\documentclass[reqno,10pt]{amsart}

\usepackage{amssymb}
\usepackage[all,cmtip]{xy}
\usepackage{xcolor}

\newtheorem{MainTheorem}{Theorem}
\newtheorem{Proposition}{Proposition}[section]
\newtheorem{Definition}[Proposition]{Definition}
\newtheorem{Lemma}[Proposition]{Lemma}
\newtheorem{Theorem}[Proposition]{Theorem}
\newtheorem{Corollary}[Proposition]{Corollary}

\DeclareMathOperator{\Val}{Val}
\DeclareMathOperator{\TVal}{TVal}
\DeclareMathOperator{\Area}{Area}
\DeclareMathOperator{\vol}{vol}
\DeclareMathOperator{\Sym}{Sym}
\DeclareMathOperator{\hol}{hol}
\DeclareMathOperator{\anti}{anti}
\DeclareMathOperator{\sgn}{sgn}
\DeclareMathOperator{\Kl}{Kl}
\DeclareMathOperator{\Gr}{Gr}
\DeclareMathOperator{\grad}{grad}
\DeclareMathOperator{\id}{id}
\DeclareMathOperator{\pd}{pd}
\DeclareMathOperator{\Dens}{Dens}
\DeclareMathOperator{\contr}{contr}
\DeclareMathOperator{\tr}{tr}
\DeclareMathOperator{\nc}{nc}
\DeclareMathOperator{\SO}{SO}
\DeclareMathOperator{\On}{O}
\DeclareMathOperator{\SL}{SL}
\DeclareMathOperator{\GL}{GL}

\DeclareMathOperator{\AGr}{\overline{Gr}}
\newcommand{\R}{\mathbb{R}}
\renewcommand{\i}{\mathbf{i}}
\newcommand\flag[2]{\left[\begin{array}{c} #1\\ #2 \end{array}
  \right]} 

\title{Kinematic formulas for tensor valuations}
\author{Andreas Bernig}
\author{Daniel Hug}
 
\email{bernig@math.uni-frankfurt.de}
\email{daniel.hug@kit.edu}

\address{Institut f\"ur Mathematik, Goethe-Universit\"at Frankfurt,
Robert-Mayer-Strasse 10, 60054 Frankfurt, Germany}
\address{Karlsruhe Institute of Technology, Department of Mathematics,  D-76128 Karls\-ruhe, Germany}

\thanks{A. B. was supported by DFG grants BE 2484/3-1 and BE 2484/5-1. D. H. was supported by DFG grants HU 1874/2-1 and HU 1847/3-1 \\ AMS 2010 {\it Mathematics subject
classification}: 52A22, 
53C65. 
}
\begin{document}

\begin{abstract}
We prove new kinematic formulas for tensor valuations and simplify previously known Crofton formulas by using the
recently developed algebraic theory of translation invariant valuations. The heart of the paper is the
computation of the Alesker-Fourier transform on the large class of spherical valuations, which is achieved by  differential-geometric and
representation theoretical tools. We also describe in explicit form the product and convolution of tensor valuations.  
\end{abstract}

\maketitle 
\tableofcontents
\section{Introduction}

Let $V$ be a finite-dimensional (real) vector space and $\mathcal{K}(V)$ the space of all compact convex subsets in $V$. Let $A$
be an abelian semigroup. A map $\phi:\mathcal{K}(V) \to A$ is called a {\it valuation} if 
\begin{displaymath}
 \phi(K \cup L)+\phi(K \cap L)=\phi(K)+\phi(L)
\end{displaymath}
whenever $K, L, K \cup L \in \mathcal{K}(V)$. If $A$ is a topological semigroup and $V$ is endowed with some Euclidean scalar product, then $\phi$ is called  continuous if it
is continuous with respect to the topology induced by the Hausdorff metric. 

Scalar valued valuations (where $A=\R$ or $A=\mathbb{C}$) have played a prominent role in integral geometry since Hadwiger's
characterization of intrinsic volumes. In particular, the kinematic formulas due to Blaschke-Chern-Santal\'o can easily be
 derived from Hadwiger's theorem \cite{klain_rota}. Recently discovered algebraic structures on the space of
translation invariant
valuations have been used to prove kinematic formulas in hermitian spaces \cite{alesker03_un, bernig_fu06, bernig_fu_hig,
bernig_fu_solanes}. 

A quite active area of research is that of {\it Minkowski valuations}, where $A=\mathcal{K}(V)$ together with the
Minkowski addition. We refer
to \cite{abardia12, abardia_bernig, haberl10, ludwig_2005, schuster10, schuster_wannerer} to cite just a few papers in
this direction. Other important cases of semigroups $A$ are the space of signed measures on $V$ (see
\cite{bernig_fu_solanes, schneider_book14} for {\it curvature measures} and corresponding {\it local kinematic
formulas}) and the space of signed measures on the unit sphere (see \cite{schneider_book14, wannerer_area_measures,
wannerer_unitary_module} for  {\it area measures} and the corresponding {\it local additive kinematic formulas}).  

In this paper, we investigate tensor valuations, where  $A=\Sym^* V$, the space of 
symmetric tensors over the given vector space $V$. Tensor
valuations were first studied by Hadwiger-Schneider \cite{hadwiger_schneider71} who considered vector valued valuations
on the basis of 
the accompanying work  
\cite{schneider71, schneider72}.
 McMullen \cite{mcmullen97} conjectured that certain tensor
valuations $Q^i\Phi_{k,r,s}$ span the space of all isometry covariant tensor valuations and proved some linear
dependences for these valuations. McMullen's conjecture was confirmed by Alesker \cite{alesker99}, 
who made essential use of his results
from \cite{alesker99_annals}. The question of determining all linear dependences among these basic tensor valuations
was settled by the
second-named author, Schneider and R.~Schuster \cite{hug_schneider_schuster_a}. 
In the subsequent paper \cite{hug_schneider_schuster_b}, the same authors studied integral-geometric formulas (Crofton-type
formulas) which are satisfied by the basic tensor valuations. Very recently, local tensor valuations were classified in 
\cite{hug_schneider_localtensor, schneider13}.

More specifically, in the present paper we  consider translation invariant tensor valuations (although we believe that a similar
study in the non-translation invariant case will be fruitful as well). In Section \ref{sec_algstruc} we generalize
 some of the  algebraic structures such as convolution, product and Alesker-Fourier transform from scalar valued valuations to
tensor valuations and establish links to new kinematic formulas. While it is not difficult to introduce these structures, 
the crucial aim  of the present work  is to determine the algebraic structures and related integral-geometric formulas
as explicitly as possible. The main technical part is the computation of the Alesker-Fourier transform of a tensor valuation.
Here we need the recently introduced Alesker-Fourier transform on odd (scalar valued) valuations \cite{alesker_fourier}. Using the 
formalism developed in \cite{bernig_fu06}, we 
provide explicit formulas for the convolution of tensor valuations.  Based on the convolution and the Alesker-Fourier transform,  
we then obtain formulas for the product of tensor
valuations as well. In the last section, we translate these results into integral-geometric formulas, that is, 
 kinematic formulas involving addition or intersection of convex sets.  

To describe our main results, we introduce some notation. Let $V$ be an $n$-dimensional Euclidean
vector space. Throughout the paper we will assume that $n \geq 2$. The space $\mathcal{K}(V)$ is
endowed with the Hausdorff metric. The set of all continuous, translation invariant valuations on $\mathcal{K}(V)$ is
denoted by $\Val(V)$, the subspace of {\it smooth} elements by $\Val^{sm}(V)$ (see Section 2 for details). If $\mu \in \Val(V)$ 
and $k\in\mathbb{N}_0$ are such that $\mu(tK)=t^k 
\mu(K)$  for all $K \in \mathcal{K}(V)$ and all $t>0$, then $\mu$ is said to be {\em homogeneous of degree $k$}. By a result of
McMullen, $\Val(V)$ can be decomposed into its homogeneous components, that is,
$$
\Val(V)=\bigoplus_{k=0}^n \Val_k(V),
$$
where $\Val_k(V)$ is the space of $k$-homogeneous elements in $\Val(V)$.

For $K \in \mathcal{K}(V)$ and $k\in\{0,\ldots,n-1\}$, let $S_k(K,\cdot)$ denote the $k$-th surface area measure of $K$, 
normalized as in \cite{schneider_book14}. In the following, we will often use that the centroid of $S_k(K,\cdot)$ is the origin, that is,
\begin{equation} \label{eq_centroid}
 \int_{S^{n-1}} y \ dS_k(K,y)=0,
\end{equation}
see \cite[(5.30)]{schneider_book14}.

The volume of the unit ball $B^n$ in $\R^n$ is denoted
by $\kappa_n:={\pi^\frac{n}{2}}/{\Gamma\left(\frac{n}{2}+1\right)}$ and the volume of the unit sphere $S^{n-1}$ in
$\R^n$ is denoted by $\omega_n=n\kappa_n$. We also define the {\it flag coefficients} 
\begin{displaymath}
 \flag{n}{k}:=\binom{n}{k} \frac{\kappa_n}{\kappa_k \kappa_{n-k}},\qquad k \in \{0,\ldots,n\}.
\end{displaymath}
The Grassmann manifold of all $k$-dimensional subspaces in a vector space $V$ is denoted by $\Gr_k(V)$. The affine
Grassmann manifold, denoted by $\AGr_k(V)$, is the space of all affine $k$-dimensional subspaces (which are also called 
$k$-flats). In the case
$V=\R^n$, $\Gr_k(V)$ is endowed with the rotation invariant probability measure. The rigid motion invariant measure on
the (non-compact)
affine Grassmann manifold $\AGr_k(V)$ is normalized in such a way that the measure of all  $k$-flats intersecting the unit
ball is  $\kappa_{n-k}$ (cf.\ \cite{schneider_book14, schneider_weil08}). 

A {\it spherical harmonic of degree $s$} in $\R^n$ is the restriction to the unit sphere of a homogeneous harmonic polynomial of
degree $s \in \mathbb{N}_0$. The space of spherical harmonics of degree $s$ in $\R^n$  is denoted by $\mathcal{H}_s^n$. 

\begin{Definition}
Let $f$ be a smooth function on $S^{n-1}$. For $k\in\{0,\ldots,n-1\}$, we define a valuation $\mu_{k,f} \in \Val_k(\R^n)$ by 
\begin{displaymath}
 \mu_{k,f}(K):=\binom{n-1}{k}\frac{1}{\omega_{n-k}}  \int_{S^{n-1}} f(y) \,dS_k(K,y).
\end{displaymath} 
\end{Definition}

If $f \in \mathcal{H}_1^n$, i.e., if $f$ is the restriction of a linear function, then $\mu_{k,f}=0$ by \eqref{eq_centroid}. 
Hence we will assume in the following that $f \perp \mathcal{H}_1^n$, that is, 
\begin{displaymath}
 \int_{S^{n-1}} f(y)y\,d\vol_{S^{n-1}}(y)=0. 
\end{displaymath}

Note that $\mu_{k,1}=\mu_k$ is the $k$-th intrinsic volume (which is also denoted by $V_k$ by some authors). Valuations
of the form $\mu_{k,f}$ are smooth and will be called {\it spherical}, since in a sense which will be made 
precise below (see Corollary \ref{Corhns}), they
correspond to spherical representations of the (proper) rotation group $\SO(n)$. 

The Alesker-Fourier transform $\mathbb{F}$ on smooth translation invariant valuations was introduced by Alesker in
\cite{alesker03_un} and \cite{alesker_fourier}. Its construction is very involved. Previously, explicit formulas were
known in the even case only. In our first main theorem, we compute the Alesker-Fourier transform of spherical valuations.

\begin{MainTheorem} \label{mainthm_fourier}
Let $f \in \mathcal{H}^n_s$, $s \neq 1$ and $1 \leq k \leq n-1$. Then 
\begin{equation} \label{eq_fourier_constant}
 \mathbb{F} (\mu_{k,f})=
\i^s
\frac{\Gamma\left(\frac{n-k}{2}\right)\Gamma\left(\frac{s+k}{2}\right)}{\Gamma\left(\frac{k}{2}\right)\Gamma\left(\frac{
s+n-k}{2}\right)} \,
\mu_{n-k,f},
\end{equation}
where $\i=\sqrt{-1}$.
\end{MainTheorem}

We denote by $\TVal^s (V)$ the vector space of translation invariant continuous valuations on ${\mathcal K}(V)$ taking their values in the vector space $\Sym^s V$ of symmetric tensors of rank $s$ over $V$. In the following, we identify $\TVal^s (V)$ with the tensor product $\Val(V)\otimes\Sym^sV$. 

Let $\TVal^{s,sm}(V)$ denote the set of smooth elements of $\TVal^s (V)$ (see Section \ref{sec_background}  for further details), which 
can be identified with $\Val^{sm}(V) \otimes \Sym^s $. If $G$ is a subgroup of the orthogonal group of
$V$, we denote by
$\TVal^{s,G}(V)$ the space of translation invariant and $G$-covariant (see Section \ref{sec_algstruc}) tensor valuations of rank $s$. 
If $G$ is a compact 
subgroup of the orthogonal group of $V$ acting 
transitively on the unit sphere, then $\TVal^{s,G}(V)$ consists of smooth valuations only and is finite-dimensional. 

We define the basic translation invariant and $O(n)$-covariant continuous tensor valuations $\Phi_{k,s} \in \TVal^{s,\On(n)}(\R^n)$, for $s\in\mathbb{N}_0$, by 
\begin{align*}
 \Phi_{k,s}(K) & := \binom{n-1}{k} \frac{1}{\omega_{n-k+s} s!} \int_{S^{n-1}} y^s \, dS_k(K,y), \quad k\in\{0,\ldots,n-1\},\\
 \Phi_{n,0}(K) & := \vol_n(K),
\end{align*}
where $y^s$ denotes the $s$-fold tensor product $y\otimes\cdots\otimes y$ of $y\in\R^n$ for $s\in\mathbb{N}$ and $y^0:=1$. Moreover, we put 
$\Phi_{n,s}:=0$ for $s\neq 0$. 
Clearly, each coefficient of $\Phi_{k,s}$ (with respect to some basis of $\Sym^s \R^n$) is a spherical valuation.  Furthermore, $\Phi_{k,0}=\mu_k$ for all $k$.

We set $Q:=\sum_{i=1}^n e_i^2 \in \Sym^2 V$, where $e_1,\ldots,e_n$ is an
orthonormal basis of $V$. Then, for each fixed $k\in\{1,\ldots,n-1\}$, the valuations $Q^i \Phi_{k,s}$ form a basis of the space of continuous, translation
invariant and $\On(n)$-covariant tensor valuations in $V=\R^n$ which are homogeneous of degree $k$ (if $n \geq 3$, then $\SO(n)$-covariance would be enough). See   \cite[Lemma 4.8 and Proposition 4.9]{alesker99_annals} for more information. 

Algebraic structures on scalar valued valuations (such as convolution, product and Alesker-Fourier transform) play an important role in integral geometry, in particular in connection with intersectional and additive kinematic formulas \cite{bernig_aig10, bernig_fu06, bernig_fu_hig, bernig_fu_solanes, fu_barcelona, wannerer_area_measures, wannerer_unitary_module}. These algebraic structures will be extended to the tensor valued case in Section \ref{sec_algstruc}. The explicit knowledge of these operations given in our second main theorem below, which summarizes Theorems \ref{thm_convolution}, \ref{thm_fourier_tensorcase} and \ref{thm_product_tensorcase}, will be the basis for writing down kinematic formulas for tensor valuations.    

Before stating Theorem \ref{mainthm_algebra_tensorvals}, we remark that for integers $a,b \geq 0$
\begin{align*}
(Q^a \Phi_{k,s_1}) \cdot (Q^b \Phi_{k,s_2}) & =Q^{a+b} (\Phi_{k,s_1} \cdot \Phi_{k,s_2}), \\
(Q^a \Phi_{k,s_1}) * (Q^b \Phi_{k,s_2})& =Q^{a+b} (\Phi_{k,s_1} * \Phi_{k,s_2}),\\
\mathbb{F}(Q^a \Phi_{k,s}) & =Q^a \mathbb{F}(\Phi_{k,s}). 
\end{align*}
Hence it is enough to write down formulas for the $\Phi_{k,s}$.
 
\begin{MainTheorem} \label{mainthm_algebra_tensorvals}
Convolution, Alesker-Fourier transform and product of the basic translation invariant, $\mathrm{O}(n)$-covariant tensor
valuations $\Phi_{k,s}$ are given by 
\begin{align*}
\Phi_{k,s_1} * \Phi_{l,s_2} & = \frac{\omega_{s_1+s_2+2n-k-l}}{\omega_{s_1+n-k}\omega_{s_2+n-l}}
\frac{(n-k)(n-l)}{2n-k-l} \cdot \\
& \qquad \cdot \binom{2n-k-l}{n-k}  \binom{s_1+s_2}{s_1} \frac{(s_1-1)(s_2-1)}{1-s_1-s_2} 
\Phi_{k+l-n,s_1+s_2}\\
\intertext{for $k,l \leq n$ with $k+l\geq n$ and $s_1,s_2\neq 1$,}
 \mathbb{F} (\Phi_{k,s}) &=\i^s \sum_{j=0}^{\lfloor \frac{s}{2}\rfloor}
\frac{(-1)^j}{(4\pi)^j j!}  Q^j
\Phi_{n-k,s-2j}\\
\intertext{for $ 0 \leq k \leq n$ and $ s \neq 1$, and}
 \Phi_{k,s_1} \cdot \Phi_{l,s_2} & = \frac{kl}{k+l} \binom{k+l}{k}  \sum_{\substack{a=0\\2a \neq s_1+s_2-1}}^{\lfloor
\frac{s_1+s_2}{2}\rfloor}
\frac{1}{(4\pi)^a a!} \Bigg( \sum_{m=0}^a \sum_{i=\max \left\{0, m-\left\lfloor\frac{s_2}{2}\right\rfloor\right\}}^{\min
\left\{m,\left\lfloor \frac{s_1}{2}\right\rfloor\right\}} \\
& \quad  (-1)^{a-m}
 \binom{a}{m} \binom{m}{i} 
\frac{\omega_{s_1+s_2-2m+k+l}}{\omega_{s_1-2i+k}\omega_{s_2-2m+2i+l}}
 \cdot \\
& \quad  \hspace{-1cm} \cdot \binom{s_1+s_2-2m}{s_1-2i} \frac{(s_1-2i-1)(s_2-2m+2i-1)}{1-s_1-s_2+2m} \Bigg) Q^a
\Phi_{k+l,s_1+s_2-2a}\allowdisplaybreaks\\
\end{align*}
for $k+l\le n$ and $s_1,s_2\neq 1$.
\end{MainTheorem}

In the extremal cases $k,l\in\{0,n\}$, some of the expressions on the right-hand side of the preceding equations have to be
interpreted properly as limits and are formally included then. For instance, if $k=n$ and $s_1 \neq 0$, then both sides of the 
first equation are zero. However, if $k=n$ and $s_1= 0$, then $(n-k)/\omega_{s_1+n-k}$ is taken as $1/\kappa_{n-k}=1/\kappa_0=1$. 
Moreover, the factors $\frac{kl}{k+l} \binom{k+l}{k}$ and $\frac{\omega_{s_1+s_2-2m+k+l}}{\omega_{s_1-2i+k}\omega_{s_2-2m+2i+l}}$ in 
the last formula may be $0$ and $\infty$ if $k=0$ or $l=0$ and will cancel out in the appropriate way. Similar remarks apply in the following.   

We note that the special case $s=s_1=s_2=0$ reduces to the following well-known  formulas for intrinsic volumes:
\begin{displaymath}
 \mathbb F \mu_k=\mu_{n-k}, \quad \mu_k \cdot \mu_l=\flag{k+l}{k}\mu_{k+l}, \quad \mu_k * \mu_l=\flag{2n-k-l}{n-k}\mu_{k+l-n},
\end{displaymath}
where $k,l$ are subject to the same conditions as in Theorem \ref{mainthm_algebra_tensorvals}.

As a first consequence,  we obtain Crofton-type formulas in the translation invariant case. More general formulas
of this type, but with rather complicated coefficients, were proved in \cite{hug_schneider_schuster_b}.

\begin{MainTheorem}[Crofton formula in the $\Phi$-basis] \label{mainthm_crofton}
If  $k,l \geq 0$ with $k+l \leq n$ and $s\in\mathbb{N}_0$, then 
\begin{align}
 \int_{\AGr_{n-l}(\R^n)} \Phi_{k,s}(K \cap \bar E)\, d\bar E & = \flag{n}{l}^{-1} \binom{k+l}{k} \frac{kl}{2(k+l)}
\frac{1}{\Gamma\left(\frac{k+l+s}{2}\right)} \cdot \nonumber\\
& \qquad \cdot \sum_{j=0}^{\lfloor
\frac{s}{2}\rfloor} \frac{\Gamma\left(\frac{l}{2}+j\right)\Gamma\left(\frac{k+s}{2}-j\right)}{(4\pi)^j\label{Croftonformula}
j!} Q^j \Phi_{k+l,s-2j}(K).
\end{align}
\end{MainTheorem} 

Note that the symmetric tensor product with $Q^a$, $a>0$, can be interchanged with the integration on the left-hand side of equation \eqref{Croftonformula} in Theorem 
\ref{mainthm_crofton}.  Moreover, if $k=0$ and $s$ is odd, then both sides of this equation are zero. If $k=0$ and $s$ is even, 
then we consider $\tfrac{k}{2}\Gamma\left(\tfrac{k}{2}\right)$ as $\Gamma\left(\tfrac{k}{2}+1\right)=1$. In case  $l=0$,  a similar discussion applies. 
Properly interpreted (and using $\frac{k}{k+l}+\frac{l}{k+l}=1$), the equation even holds for $k=l=0$. 


\medskip

In Section \ref{sec_algstruc} we will extend the {\it fundamental theorem of algebraic integral geometry}
\cite{bernig_aig10, fu_barcelona} to tensor valuations. Roughly speaking, it explains how the algebraic
structures from our second main theorem determine the integral-geometric formulas. As an application, we compute the intersectional kinematic formulas
for some tensors of small rank. However, our approach makes clear how such formulas can be determined for tensors of arbitrary rank.  

Let $\overline{\On(n)} :=\On(n) \ltimes \R^n$ be the Euclidean motion group, endowed with the product measure (as
in \cite{schneider_weil08}). For $\bar g \in \On(n)$, we let $g \in \mathrm{O}(n)$ be the rotation part. 
Then the intersectional kinematic operator 
\begin{displaymath}
 k_{s_1,s_2}^{\On(n)}:\TVal^{s_1+s_2,\On(n)}(\R^n) \to \TVal^{s_1,\On(n)}(\R^n) \otimes
\TVal^{s_2,\On(n)}(\R^n)
\end{displaymath}
is defined by 
\begin{displaymath}
 k_{s_1,s_2}^{\On(n)}(\Phi)(K,L):=\int_{\overline{\On(n)}} (\id^{s_1} \otimes g^{s_2})\Phi(K \cap  \bar g^{-1}
L)\, d\bar g. 
\end{displaymath}
Here $\id^{s_1} \otimes g^{s_2}$ acts in the natural way on $\Sym^{s_1+s_2}(\R^n) \subset (\R^n)^{\otimes s_1}
\otimes (\R^n)^{\otimes s_2}$.

\begin{MainTheorem} \label{mainthm_kf}
For $0 \leq i \leq
n-1$, the intersectional kinematic formulas of bi-rank $(2,2)$, $(3,2)$ and $(3,3)$ respectively are given by 
\begin{align*}
k_{2,2}^{\On(n)}(\Phi_{i,4})  & = \frac{1}{48 \pi^2 \Gamma\left(\frac{n+3}{2}\right)\Gamma\left(\frac{i+1}{2}\right)}
\sum_{k+l=n+i} \frac{\Gamma\left(\frac{k+1}{2}\right)\Gamma\left(\frac{l+1}{2}\right)}{
kl } \cdot\\
& \quad \cdot\bigg( 4\pi^2 (ni^2+i^2-2ni-2i-2n) \Phi_{k,2} \otimes \Phi_{l,2}\\
& \quad\qquad -\pi (ik+nik-2k-3ni-n^2i) \Phi_{k,2}
\otimes Q\Phi_{l,0}\\
& \quad\qquad -\pi (il+nil-2l-3ni-n^2i) Q\Phi_{k,0} \otimes \Phi_{l,2}\\
& \quad\qquad +\frac14 (n+3)(i-l)(i-k) Q\Phi_{k,0}
\otimes Q\Phi_{l,0}\bigg),\\
 k_{3,2}^{\On(n)}(\Phi_{i,5})  & =  \sum_{k+l=n+i}
\frac{(i+1)\Gamma\left(\frac{l+1}{2}\right)\Gamma\left(\frac{k}{2}\right)}{40 \pi
(k+1)l \Gamma\left(\frac{n+1}{2}\right) \Gamma\left(\frac{i}{2}\right)}\cdot \\
&  \qquad \cdot \Phi_{k,3} \otimes \left(4 \pi (i-3) \Phi_{l,2}+(n-k+3) Q\Phi_{l,0}\right),\\
 k_{3,3}^{\On(n)}(\Phi_{i,6}) & =\frac{(i+1)(i-1)(i-3)}{40 \Gamma\left(\frac{n+1}{2}\right)\Gamma\left(\frac{i+1}{2}\right)}
\sum_{k+l=n+i} \frac{\Gamma\left(\frac{k}{2}\right)\Gamma\left(\frac{l}{2}\right)
}{(k+1)(l+1)} \Phi_{k,3} \otimes \Phi_{l,3}.
\end{align*}
\end{MainTheorem}

The kinematic formulas for the group $\SO(n)$ with $ n \geq 3$ are formally identical, since
$\TVal^{\On(n)}(\R^n)
\cong \TVal^{\SO(n)}(\R^n)$ for $n \geq 3$ (see \cite{alesker99_annals}).

We also point out that the strategy of proof for this theorem can be used to derive formulas for $k_{s_1+a_1,s_2+a_2}^{\On(n)}(Q^a \Phi_{i,s_1+s_2})$ with $2a=a_1+a_2>0$.

Now we turn to additive kinematic formulas for tensor valuations. In this case, we do not need the algebraic machinery developed in earlier sections, since we can relate the additive kinematic formulas for tensor valuations to the well-known additive kinematic formulas for area measures. 

To state the result, we introduce the additive kinematic operator 
\begin{displaymath}
 a_{s_1,s_2}^{\On(n)}:\TVal^{s_1+s_2,\On(n)}(\R^n) \to \TVal^{s_1,\On(n)}(\R^n) \otimes
\TVal^{s_2,\On(n)}(\R^n)
\end{displaymath}
as 
\begin{displaymath}
 a_{s_1,s_2}^{\On(n)}(\Phi)(K,L):=\int_{\On(n)} (\id^{s_1} \otimes g^{s_2})\Phi(K +  g^{-1}
L)\, d g,
\end{displaymath}
where the integration over $\On(n)$ is with respect to the Haar probability measure. 
The result is the following theorem.

\begin{MainTheorem}[Additive kinematic formula for tensor valuations] \label{thm_additive_kf}
Let $0 \leq i \leq n-1$ and $s_1,s_2 \neq 1$. Then 
\begin{align*}
a_{s_1,s_2}^{\On(n)}(\Phi_{i,s_1+s_2})
& = \binom{s_1+s_2}{s_1}^{-1} \sum_{k+l=i} \frac{(n-k-1)!(n-l-1)!}{(n-i-1)!(n-1)!}\cdot\\
&\qquad \cdot \frac{\omega_{n-k+s_1}
\omega_{n-l+s_2}}{\omega_n \omega_{n-i+s_1+s_2}} \Phi_{k,s_1} \otimes \Phi_{l,s_2}.
\end{align*}
\end{MainTheorem}

Here we do not state additive kinematic formulas for the valuations $Q^a \Phi_{i,s_1+s_2}$ with $ a>0$. We do not know whether such formulas can be deduced from  additive kinematic formulas for area measures. Nevertheless, there is enough information in Theorems \ref{mainthm_algebra_tensorvals} and \ref{thm_ftaig} to write down such formulas as well (using a similar approach as in the proof of Theorem \ref{mainthm_kf}).

\subsection*{Plan of the paper}
In Section \ref{sec_background}, we will review the necessary background from algebraic integral geometry of translation invariant 
scalar valued valuations. The
extension of this theory to translation invariant tensor valuations is the subject of Section \ref{sec_algstruc}. There we
will introduce convolution, Alesker-Fourier transform and product of translation invariant tensor valuations as well as additive
and intersectional kinematic formulas. The relation between the algebraic structures and the kinematic formulas is
explored in Theorem \ref{thm_ftaig}. This relationship is illustrated by a commutative diagram which explains how convolution and additive kinematic formulas on the one hand, and multiplication 
and intersectional kinematic formulas on the other hand, are related via Poincar\'e 
duality and Fourier transform. 

The heart of the paper is Section \ref{sec_fourier}, where we explicitly compute the Alesker-Fourier
transform of spherical valuations, which is the key step for determining the  Fourier transform of 
continuous, translation invariant and $\mathrm{O}(n)$-covariant tensor valuations. In this section, we consider
scalar valued rather than tensor valued valuations, but the main result, which is Theorem 
\ref{mainthm_fourier}, is then translated into the language of tensor
valuations in Section \ref{sec_alg_prop}. However, already in Section \ref{sec_fourier} it will be convenient  
to take advantage of the close relationship between spherical valuations and tensor valuations. Since the explicit determination of  
the Fourier transform of spherical valuations requires several steps, we start with the relatively simple case of even valuations, where 
the Klain map can be efficiently used. The two-dimensional case, which allows a direct approach, is treated next. Then we give a fairly detailed 
outline of our approach in the general case which explains how the subsequent arguments and calculations are connected (see Subsection \ref{subsec_outline}).

In Section \ref{sec_alg_prop}, we also compute convolution, product and Poincar\'e
duality for tensor valuations. The explicit form of the Poincar\'e pairing for tensor valuations 
will be crucial for the proof of intersectional kinematic formulas.

The final Section \ref{sec_explicit_kinforms} is devoted to kinematic formulas. First, we simplify in the translation invariant case the constants
in the general Crofton formula established in \cite{hug_schneider_schuster_b}.  
Then we prove additive kinematic 
formulas for tensor valuations. The intersectional kinematic formulas are more difficult, and we will state them in
explicit form for small bi-ranks only. However, it will be clear from the proof how formulas for tensors of higher
rank can be computed as well.

\subsection*{Acknowledgments}
We would like to thank Franz Schuster for useful remarks on an earlier version of this paper.  We are grateful to the referee for his remarks which lead to an improved presentation.

\section{Background from translation invariant valuations}
\label{sec_background}

The theory of scalar valued, translation invariant, continuous valuations is a very active and rich one, compare the
surveys \cite{bernig_aig10, fu_barcelona}. For the reader's convenience, we collect those facts which are relevant in
the present paper. 

The space of continuous, translation invariant valuations on an $n$-dimensional
vector space $V$ is denoted by $\Val(V)$ or just by $\Val$
if there is no risk of confusion. Similarly, we write $\mathcal{K}$ instead of $\mathcal{K}(V)$. 
A valuation $\phi \in \Val$ is called {\it even} if $\phi(-K)=\phi(K)$, for all $K\in\mathcal{K}(V)$, and {\it odd}
if $\phi(-K)=-\phi(K)$, for all $K\in\mathcal{K}(V)$. If there is some $k\in\mathbb{N}_0$ such 
that $\phi(tK)=t^k \phi(K) $ for all $t>0$ and  all $K\in\mathcal{K}(V)$, then $\phi$ is
said to be {\it homogeneous of degree $k$}. The space of even/odd valuations of degree $k$ is denoted by
$\Val_k^\pm(V)$. A
fundamental result by McMullen
\cite{mcmullen77} is the decomposition 
\begin{equation}\label{McMa}
 \Val(V) = \bigoplus_{\substack{k=0,\ldots,n\\ \epsilon=\pm}} \Val_k^\epsilon(V). 
\end{equation}
 
The space $\Val_0(V)$ is spanned by the Euler characteristic $\chi$ (which satisfies $\chi(K)=1$ for all $K\in\mathcal{K}(V)\setminus\{\emptyset\}$), while
$\Val_n(V)$ is spanned by a Lebesgue measure $\vol_n$ on $V$. 

The Euler-Verdier involution $\sigma$ is defined in \cite{alesker_val_man2} by
\begin{displaymath}
 (\sigma \phi)(K):=(-1)^k \phi(-K), \quad \phi \in \Val_k(V), K\in\mathcal{K}(V). 
\end{displaymath}

Let $V \cong \R^n$ be an $n$-dimensional Euclidean vector space. 
Let $SV$ be the sphere bundle over $V$ and $\nc(K)$ the
normal cycle 
(which is sometimes also called the normal bundle) of $K \in
\mathcal{K}(V)$ \cite{zaehle86}. A valuation $\phi \in \Val(V)$ is called {\it smooth} 
if there exist translation invariant differential
forms $\gamma \in \Omega^{n-1}(SV)$ and $ \tau \in \Omega^n(V)$ such that 
\begin{displaymath}
 \phi(K)=\int_{\nc(K)} \gamma+\int_K \tau. 
\end{displaymath}
The space $\Val^{sm}(V)$ of smooth valuations is dense in $\Val(V)$, compare \cite{alesker_val_man1}.

Alesker \cite{alesker04_product} has shown that $\Val^{sm}(V)$ has the structure of an associative, commutative graded algebra satisfying a version
of
Poincar\'e duality. The product of two smooth valuations satisfies the following formula (and, together with continuity
and bilinearity is characterized by this formula): if $\phi_i(K)=\vol_n(K+A_i)$, $i=1,2$, where $A_i$ is a compact convex
body with smooth boundary and positive curvature and $K\in\mathcal{K}(V)$, then 
$$
(\phi_1 \cdot \phi_2)(K)=\vol_{2n}(\Delta K+A_1 \times A_2), 
$$
where
$\Delta:V \to V \times V$ is the diagonal embedding. The Euler characteristic $\chi$ is the identity for Alesker's multiplication.

In \cite{bernig_fu06} a {\it convolution product} was introduced on $\Val^{sm}(V)$. Without referring to a Euclidean
structure, the convolution is a product on the space $\Val^{sm}(V) \otimes \Dens(V^*)$, where $\Dens(V^*)$ is the
$1$-dimensional space of complex Lebesgue measures on $V^*$. If $\phi_1,\phi_2$ are as above and $K\in\mathcal{K}(V)$, then 
$$
(\phi_1 *\phi_2)(K)=\vol_n(K+A_1+A_2). 
$$
Together with continuity and bilinearity, the convolution is determined by this formula. The neutral element of the convolution is the $n$-dimensional volume $\vol_n$.
 
The picture is completed by a Fourier-type transform, due to Alesker
\cite{alesker03_un, alesker_fourier}, which is a linear map $\mathbb{F}:\Val^{sm}(V) \to \Val^{sm}(V^*) \otimes
\Dens(V)$ satisfying 
\begin{equation} \label{eq_fourier_convprod_scalar}
\mathbb{F}(\phi_1 \cdot \phi_2)=\mathbb{F}(\phi_1) * \mathbb{F}(\phi_2)
\end{equation}
and with $\mathbb{F}^2:=\mathbb{F}\circ \mathbb{F}$ also the Plancherel-type
formula 
\begin{displaymath}
(\mathbb{F}^2 (\phi))(K)=\phi(-K), \quad K\in\mathcal{K}(V).
\end{displaymath}

Let $\phi \in \Val_k^+(V)$. If $E \in \Gr_k(V)$, then the restriction of $\phi$
to $E$ is a multiple of the Lebesgue measure on $E$. We thus obtain a function $\Kl_\phi \in C(\Gr_k(V))$, called the
{\it Klain
function of $\phi$}. The resulting map $\Val_k^+(V) \to C(\Gr_k(V))$ is injective, as was shown by Klain \cite{klain00}.
It is
called {\it Klain embedding}. The Alesker-Fourier transform on even valuations $\phi \in \Val_k^+(V)$ is characterized by
\begin{displaymath}
 \Kl_{\mathbb{F}( \phi)}(E)=\Kl_\phi(E^\perp),\quad E \in \Gr_{n-k}(V),
\end{displaymath}
if $V$ and $V^*$ are identified via the scalar product. 
Although there is an analogous embedding due to Schneider \cite{schneider96} in the odd case, no similar
characterization of the Alesker-Fourier transform on odd valuations is known. In a central part of the present paper (Section
\ref{sec_fourier}), we will compute explicitly the Alesker-Fourier transform of certain odd valuations.

\section{Kinematic formulas for tensor valuations}
\label{sec_algstruc}

\subsection{Algebraic structures on tensor valuations}

Let $V$ be a finite-dimensional Euclidean vector space of dimension $n \geq 2$. The space of symmetric tensors of rank
$s \in\mathbb{N}_0$ over $V$ will be denoted by $\Sym^s
V$, 
which is a subspace of $V^{\otimes s}$, the $s$-th tensor power of $V$. Here we put $V^{\otimes 0}:=\R$. Moreover, we define 
$\Sym^*V:=\bigoplus_{s=0}^\infty \Sym^s V$. A translation invariant (continuous) 
tensor valuation is a translation invariant (continuous) valuation with values in $\Sym^*V$. 

The space $\TVal(V)$ of all translation
invariant continuous tensor valuations will be identified with
\begin{displaymath} 
\Val(V) \otimes \Sym^*V. 
\end{displaymath}
Similarly, the space $\TVal^{sm}(V)$ of smooth translation invariant tensor valuations will be identified with 
\begin{displaymath}
 \Val^{sm} (V) \otimes \Sym^*V.
\end{displaymath}

We say that $\Phi\in \TVal(V)$ has rank $s$ if $\Phi$ has values in $\Sym^s V$, the corresponding space will be denoted by
$\TVal^s(V)$. The set of smooth elements of $\TVal^s(V)$ is denoted by $\TVal^{s,sm}(V)$.  
McMullen's decomposition \eqref{McMa} of $\Val(V)$ can be extended to a decomposition 
\begin{displaymath}
 \TVal^s(V) = \bigoplus_{\substack{k=0,\ldots,n\\ \epsilon=\pm}} \TVal_k^{s,\epsilon}(V). 
\end{displaymath}
 
If $G$ is a group acting linearly on $V$, then there is a natural induced action on all tensor powers and on $\Sym^*V$.
The
action of $G$ on $\TVal^s(V)$ is given by 
$$
(g\Phi)(K):=g^{\otimes s}(\Phi(g^{-1}K)),\qquad K\in\mathcal{K}(V),\, g\in G.
$$ 
A valuation $\Phi\in \TVal(V)$ is called {\it $G$-covariant} if
$g\Phi=\Phi$. The subspace of $G$-covariant tensor valuations of rank $s$ and degree $k$ will be denoted by
$\TVal_k^{s,G}(V)$. It is shown as in the scalar valued case \cite{fu_barcelona} that $\TVal^{s,G}(V)$ is finite-dimensional
and 
consists of smooth valuations provided $G$ is a compact subgroup of  $\On(n)$ acting transitively
on the unit sphere of $V$ (see \cite[Lemma 4.1]{wannerer_area_measures}). 

In this paper, we will be interested in the case $G=\On(n)$ (cf.\
\cite{alesker99_annals,alesker99}). In the literature \cite{hug_schneider_schuster_a,
hug_schneider_schuster_b}, $\On(n)$-covariant tensor valuations are sometimes just called {\it tensor valuations}, but
in the recent work \cite{wannerer_area_measures} it turned out that the study of $G$-covariant tensor valuations for other
groups (in particular, for the unitary group) is very fruitful as well. 

Convolution, product, Alesker-Fourier transform and Euler-Verdier involution on translation invariant valuations can be
extended component-wise to corresponding algebraic structures on $\TVal^{sm}(V)$.  If $\Phi_1(K)=\sum_{i=1}^m \phi_i(K)w_i$, where
$w_1,\ldots,w_m$ is a basis of $\Sym^{s_1}V$, and $\Phi_2(K)=\sum_{j=1}^l \psi_j(K)u_j$, where $u_1,\ldots,u_l$ is a
basis of $\Sym^{s_2}V$, then 
\begin{displaymath}
 \left(\Phi_1 \cdot \Phi_2\right)(K)=\sum_{i,j} (\phi_i \cdot \psi_j)(K) w_i   u_j. 
\end{displaymath}
The dot on the right-hand side is the product of smooth valuations and $w_i u_j\in \Sym^{s_1+s_2}V$ denotes the
symmetric tensor product of the symmetric tensors $w_i\in \Sym^{s_1}V$ and $u_j\in \Sym^{s_2}V$. 
Note that this definition is independent of the chosen bases. 
 In a similar way, convolution, Alesker-Fourier transform and Euler-Verdier involution for tensor valuations can be defined.

If $\Phi_1 \in \TVal_k^{s_1,sm}(V)$, $\Phi_2 \in \TVal_l^{s_2,sm}(V)$ and $k,l\in\{0,\ldots,n\}$, then 
\begin{align*}
  \Phi_1 * \Phi_2 & \in \TVal_{k+l-n}^{s_1+s_2,sm}(V),\quad k+l\ge n,\\
\Phi_1 \cdot \Phi_2 & \in \TVal_{k+l}^{s_1+s_2,sm}(V),\quad k+l\le n,\\
\mathbb{F} (\Phi_1) & \in \TVal_{n-k}^{s_1,sm}(V),\\
 \sigma (\Phi_1) & \in \TVal_{k}^{s_1,sm}(V).
\end{align*}
Furthermore, \eqref{eq_fourier_convprod_scalar} generalizes to 
\begin{equation} \label{eq_fourier_convprod}
 \mathbb{F}(\Phi_1 \cdot \Phi_2)=\mathbb{F}(\Phi_1) * \mathbb{F}(\Phi_2).
\end{equation}
Multiplication with powers of $Q$ commutes with these operations, for instance, we have $(Q^i \Phi_1) * \Phi_2=Q^i
(\Phi_1 * \Phi_2)$ for $i\in \mathbb{N}_0$. 

We let 
\begin{displaymath}
 m,c:\TVal^{s_1,sm}(V) \otimes \TVal^{s_2,sm}(V) \to \TVal^{s_1+s_2,sm}(V)
\end{displaymath}
be the maps corresponding to the product and the convolution. 

The space $\TVal^G(V)$ is closed under these operations. We will write $m_G,c_G$ for the restrictions to the invariant
subspaces. The adjoint maps will be denoted by $m_G^*,c_G^*$.

Unlike the 
case of $\On(n)$-invariant and translation invariant scalar valued valuations, {the coefficients of} $\On(n)$-covariant tensor valuations
are not necessarily even. The computation of the Alesker-Fourier transform will be much more involved in the odd case and forms
the heart of the present paper. 

\subsection{Kinematic formulas}

Let $G$ be a subgroup of $\On(n)$ which acts transitively on the unit sphere. Then the space $\TVal^G(V)$ is 
finite-dimensional. 
 
Let $\Phi \in \TVal^{s_1+s_2,G}(V)$. Following \cite{wannerer_area_measures}, we define a bivaluation with values in
$\Sym^{s_1}V \otimes \Sym^{s_2}V$ by 
\begin{displaymath}
 a^G_{s_1,s_2}(\Phi)(K,L):=\int_G (\id^{\otimes s_1} \otimes g^{\otimes s_2}) \Phi(K+g^{-1}L)\, dg
\end{displaymath}
for $K,L\in\mathcal{K}(V)$, where $G$ is endowed with the Haar probability measure. 

Let $\Phi_1,\ldots,\Phi_{m_1}$ be a basis of $\TVal^{s_1,G}(V)$ and let $\Psi_1,\ldots,\Psi_{m_2}$ be a basis of $\TVal^{s_2,G}(V)$. 
Then,  by the classical Hadwiger argument (compare \cite{klain_rota} or \cite[Theorem 4.3]{wannerer_area_measures}), there are constants $c_{ij}^\Phi$ such that 
\begin{displaymath}
 a^G_{s_1,s_2}(\Phi)(K,L)=\sum_{i,j} c_{ij}^\Phi\, \Phi_i(K) \otimes \Psi_j(L)
\end{displaymath}
for $K,L\in\mathcal{K}(V)$.
 
The  {\it additive kinematic operator} is the map
\begin{align*}
 a^G_{s_1,s_2}: \TVal^{s_1+s_2,G}(V) & \to \TVal^{s_1,G}(V) \otimes \TVal^{s_2,G}(V)\\
\Phi & \mapsto \sum_{i,j} c_{ij}^\Phi \,\Phi_i \otimes \Psi_j, 
\end{align*}
which is independent of the choice of the bases. 

Similarly, we can define a bivaluation with values in
$\Sym^{s_1}V \otimes \Sym^{s_2}V$ by 
\begin{displaymath}
 k^G_{s_1,s_2}(\Phi)(K,L):=\int_{\bar G} (\id^{\otimes s_1} \otimes g^{\otimes s_2}) \Phi(K \cap  \bar
g^{-1}L)\, d\bar g,
\end{displaymath}
for $K,L\in\mathcal{K}(V)$, where $\bar G$ is the group generated by $G$ and the translation group, endowed with the product measure, and where $g$ is the
rotational part of $\bar g$. Choosing bases and arguing as above, we obtain a linear map, called {\it intersectional kinematic operator}, 
\begin{align*}
 k^G_{s_1,s_2}: \TVal^{s_1+s_2,G}(V) & \to \TVal^{s_1,G}(V) \otimes \TVal^{s_2,G}(V)\\
 \Phi & \mapsto \sum_{i,j} d_{ij}^\Phi \,\Phi_i \otimes \Psi_j,
\end{align*}
which is independent of the choice of the bases. Note that the constants $d_{ij}^\Phi$ depend on $\Phi$ and the chosen bases. Hence,
\begin{equation}
 k^G_{s_1,s_2}(\Phi)(K,L)=\sum_{i,j} d_{ij}^\Phi\, \Phi_i(K) \otimes \Psi_j(L)\label{kinop}
\end{equation}
for $K,L\in\mathcal{K}(V)$.

\subsection{Relation between algebraic structures and kinematic formulas}

We refer to \cite{bernig_aig10,
bernig_fu06, bernig_fu_hig, bernig_fu_solanes,
fu_barcelona} for the relation between algebraic structures on scalar valued valuations and integral-geometric
formulas.

Let $V$ be a Euclidean vector space with scalar product $\langle\cdot,\cdot\rangle$. In the following, we will briefly write
 $\Val$ for $\Val(V)$ and  $\TVal$ for $\TVal(V)$. For $s \leq r$ we define the
contraction map by 
\begin{align*}
 \contr: V^{\otimes s} \times V^{\otimes r} & \to V^{\otimes (r-s)},\\
(v_1 \otimes \ldots \otimes v_s,w_1 \otimes \ldots \otimes w_r)
& \mapsto \langle v_1,w_1\rangle \langle v_2,w_2\rangle\ldots \langle v_s,w_s\rangle w_{s+1} \otimes \cdots \otimes
w_r
\end{align*}
and linearity.
This map restricts to a map $\contr:\Sym^s V \times \Sym^r V \to \Sym^{r-s}V$. In particular, if
$r=s$, the
map $\Sym^s V \times \Sym^s V \to \R$ is the usual scalar product on $\Sym^s V$, which will also
be denoted by $\langle
\cdot,\cdot \rangle$. For symmetric tensors $u,v$, we denote the symmetric tensor product by $uv$. 
We will make  use of the equation
\begin{displaymath}
\contr(u,\contr(v,w))=\contr(uv,w),  
\end{displaymath}
where $u \in \Sym^s V$, $v \in \Sym^r V$, $w \in \Sym^t V$ and $t \geq s+r$ \cite{wannerer_area_measures}.

The trace map $\tr:\Sym^s V \to \Sym^{s-2} V$ is defined by restriction of the map $V^{\otimes s}\to V^{\otimes (s-2)}$, 
$ v_1 \otimes \ldots \otimes v_s \mapsto \langle v_1,v_2\rangle v_3 \otimes
\ldots \otimes v_s$, for $s\ge 2$. It is obtained by contraction with $Q$, that is, $\tr w=\contr(Q,w)$ (recall
that $Q=\sum_i e_i^2$ is the metric tensor, where $e_1,\ldots,e_n$ is an orthonormal basis of $V$). We thus get 
\begin{equation} \label{eq_adjoint_Q}
 \langle Qv,w\rangle=\langle v, \tr w\rangle,\quad v \in \Sym^s V, w \in \Sym^{s+2}V.
\end{equation}

By Alesker's Poincar\'e duality \cite{alesker04_product}, there are isomorphisms $\pd_c, \pd_m :\Val^{sm} \to
(\Val^{sm})^*$ defined by 
\begin{align*}
 \langle \pd_c
\phi,\psi\rangle & :=(\phi * \psi)_0, \quad \phi,\psi \in \Val^{sm},\\
\langle \pd_m
\phi,\psi\rangle & :=(\phi \cdot \psi)_n, \quad \phi,\psi \in \Val^{sm}.
\end{align*}
Here the lower index $0$ stands for the coefficient of $\chi$ in the decomposition into homogeneous components, while
the lower index $n$ stands for the coefficient of $\vol_n$.  

The scalar product on $\Sym^s V$ induces an isomorphism $q^s:\Sym^s V \to (\Sym^s V)^*$ and we set 
\begin{align*}
 \pd^s_c & : \TVal^{s,sm}=\Val^{sm} \otimes \Sym^s V
 \xrightarrow[]{\pd_c \otimes q^s}
(\Val^{sm})^* \otimes (\Sym^s
V)^*=(\TVal^{s,sm})^*,\\
 \pd^s_m & : \TVal^{s,sm}=\Val^{sm} \otimes \Sym^s V 
 \xrightarrow[]{\pd_m \otimes q^s} (\Val^{sm})^*
\otimes (\Sym^s
V)^*=(\TVal^{s,sm})^*.
\end{align*}

\begin{Lemma} \label{lemma_relation_pds}
\begin{displaymath}
\pd_m^s=(-1)^s \pd_c^s.
\end{displaymath}
\end{Lemma}

\proof  
It was shown in  \cite[(11), (12)]{wannerer_area_measures} that for smooth translation invariant valuations $\phi,\psi$ on $\R^n$ 
\begin{displaymath}
 \langle \pd_m(\phi),\psi\rangle=\left\{ \begin{matrix} \langle \pd_c(\phi),\psi\rangle, & \phi,\psi \text{ even},\\ -\langle \pd_c(\phi),\psi\rangle, & \phi,\psi \text{ odd.}\end{matrix}\right.
\end{displaymath}
If $\phi,\psi$ have opposite parity, then both pairings vanish. Hence the statement follows from the observation that each coefficient of a tensor valuation of rank $s$ is of parity $(-1)^s$. 
\endproof

\begin{Theorem} \label{thm_ftaig}
Let $G$ be a compact subgroup of $\On(n)$ acting transitively on the unit sphere. Then the  diagram
\begin{displaymath}
 \xymatrix{\TVal^{s_1+s_2,G} \ar[r]^<<<<<<<<<{a^G_{s_1,s_2}} \ar[d]_{\pd_c^{s_1+s_2}} & \TVal^{s_1,G} \otimes
\TVal^{s_2,G}
\ar[d]_{\pd_c^{s_1} \otimes \pd_c^{s_2}} \\
\left(\TVal^{s_1+s_2,G}\right)^* \ar[r]^<<<<<{c_G^*} \ar[d]_{\mathbb{F}^*} & \left(\TVal^{s_1,G}\right)^* \otimes
\left(\TVal^{s_2,G}\right)^* \ar[d]_{\mathbb{F}^* \otimes \mathbb{F}^*}\\
\left(\TVal^{s_1+s_2,G}\right)^* \ar[r]^<<<<<{m_G^*} & \left(\TVal^{s_1,G}\right)^* \otimes
\left(\TVal^{s_2,G}\right)^*\\
\TVal^{s_1+s_2,G} \ar[r]^<<<<<<<<<{k^G_{s_1,s_2}} \ar[u]^{\pd_m^{s_1+s_2}} & \TVal^{s_1,G}
\otimes \TVal^{s_2,G}
\ar[u]^{\pd_m^{s_1} \otimes \pd_m^{s_2}}}
\end{displaymath}
commutes and encodes the relations between product, convolution,
Alesker-Fourier transform, intersectional and additive kinematic formulas.
\end{Theorem}

\proof
The commutativity of the upper square was shown in \cite[Theorem 4.6]{wannerer_area_measures}. The commutativity of the
middle square is equivalent to
\eqref{eq_fourier_convprod}. 

Let us show the commutativity of the lower square. 

For this, let $\kappa \in \Sym^{s_1}V$, $\lambda \in \Sym^{s_2}V$ and let $A,B \in \mathcal{K}(V)$ 
be smooth with positive curvature. We
define valuations $\mu_{A,\kappa}:=\vol(\bullet-A) \otimes \kappa \in \TVal^{s_1,sm}$ and
$\mu_{B,\lambda}:=\vol(\bullet-B) \otimes \lambda \in \TVal^{s_2,sm}$.

Looking at the components with respect to bases of $\Sym^{s_1}V$ and $ \Sym^{s_2}V$, one easily checks that 
\begin{align*}
 \langle \pd_m^{s_1} \Phi,\mu_{A,\kappa}\rangle & = \langle \Phi(A),\kappa\rangle, \quad \Phi \in \TVal^{s_1,sm},\\
 \langle \pd_m^{s_2} \Phi,\mu_{B,\lambda}\rangle & = \langle \Phi(B),\lambda\rangle, \quad \Phi \in \TVal^{s_2,sm},\\
\left\langle \pd_m^{s_1+s_2} \Phi,\mu_{A,\kappa} \cdot \mu_{B,\lambda}\right\rangle & = \left\langle \int_V \Phi((y+A)
\cap B)\, dy, \kappa  \lambda\right\rangle, \quad \Phi \in \TVal^{s_1+s_2,sm}.
\end{align*}

Define 
\begin{align*}
\tilde \mu_{A,\kappa}&:=\int_G \mu_{gA,g\kappa} \, dg \in \TVal^{s_1,G},\\
\tilde \mu_{B,\lambda}&:=\int_G \mu_{gB,g\lambda} \, dg \in \TVal^{s_2,G}.
\end{align*} 
Alesker's irreducibility theorem \cite{alesker_mcullenconj01}
implies that linear
combinations of valuations of the form $\tilde \mu_{A,\kappa}$ span $\TVal^{s_1,G}$ and similarly for $\TVal^{s_2,G}$
(compare also \cite{fu_barcelona} for this argument).

By the continuity of the product, we deduce that
\begin{align}
 \langle \pd_m^{s_1} \Phi,\tilde \mu_{A,\kappa}\rangle & = \int_{G}\langle \Phi(gA),g\kappa\rangle\, dg, \quad \Phi \in \TVal^{s_1,sm},\label{(a)}\\
 \langle \pd_m^{s_2} \Phi,\tilde \mu_{B,\lambda}\rangle & = \int_{G}\langle \Phi(hB),h\lambda\rangle\, dh, \quad \Phi \in \TVal^{s_2,sm},\label{(b)}\\
\left\langle \pd_m^{s_1+s_2} \Phi,\tilde\mu_{A,\kappa} \cdot \tilde\mu_{B,\lambda}\right\rangle & = \int_{G}\int_{G}\left\langle \int_V \Phi((y+gA)
\cap hB)\, dy, (g\kappa)  (h\lambda)\right\rangle\,dg\, dh, \label{(c)}
\end{align}
for $\Phi \in \TVal^{s_1+s_2,sm}$.
 
Let $\Phi \in \TVal^{s_1+s_2,G}$. Using \eqref{kinop}, \eqref{(a)}, \eqref{(b)}, and then again \eqref{kinop}, we obtain
\begin{align*}
 \langle (\pd_m^{s_1} \otimes & \pd_m^{s_2}) \circ k^G_{s_1,s_2}(\Phi), \tilde \mu_{A,\kappa} \otimes
\tilde \mu_{B,\lambda}\rangle \\
& = \sum_{i,j}d_{ij}^{\Phi}\langle \pd_m^{s_1}\Phi_i,\tilde\mu_{A,\kappa}\rangle \langle \pd_m^{s_2}\Psi_j,\tilde\mu_{B,\lambda}\rangle \\
& = \sum_{i,j}d_{ij}^{\Phi}\int_G\langle \Phi_i(gA),g\kappa\rangle \, dg \int_G\langle \Psi_j(hB),h\lambda\rangle \, dh\\
& = \int_G \int_G \langle k^G_{s_1,s_2}(\Phi)(gA,hB),g\kappa \otimes h\lambda\rangle \, dg\, dh.
\end{align*}
Next we use the definition of the kinematic operator, the covariance of $\Phi$, the invariance properties of the integrating measure on $\overline{G}$ and the 
normalization of the invariant measure on $G$ to arrive at
\begin{align*}
\langle (\pd_m^{s_1} \otimes & \pd_m^{s_2}) \circ k^G_{s_1,s_2}(\Phi), \tilde \mu_{A,\kappa} \otimes
\tilde \mu_{B,\lambda}\rangle \\
& = \int_G \int_G \int_{\overline{G}} \left\langle (\mathrm{id}^{s_1} \otimes r^{s_2}) \Phi(gA \cap \bar r^{-1}
hB),g\kappa  \otimes
h\lambda \right\rangle \, d\bar r\, dg \, dh \\
& = \int_G \int_G \int_{\overline{G}} \left\langle \Phi(gA \cap \bar r^{-1} hB),g\kappa  \otimes
(r^{-1} \circ h)\lambda \right\rangle \,d\bar r\, dg\, dh \\
& = \int_{\overline{G}} \left\langle \Phi(A \cap \bar g^{-1} B),\kappa \otimes
g^{-1}\lambda\right\rangle\, d\bar g\\
& = \int_{\overline{G}} \left\langle \Phi(A \cap \bar g^{-1} B),\kappa 
g^{-1}\lambda\right\rangle\, d\bar g.
\end{align*}

On the other hand,  by \eqref{(c)} and since $\overline{G} = G \ltimes V$ is endowed with the product measure,
\begin{align*}
 \langle m_G^* \circ &\pd_m^{s_1+s_2}(\Phi) ,\tilde \mu_{A,\kappa} \otimes \tilde \mu_{B,\lambda}\rangle\\
& =\langle
\pd_m^{s_1+s_2}(\Phi),\tilde \mu_{A,\kappa} \cdot \tilde \mu_{B,\lambda}\rangle\\
%
%
& = \int_G \int_G \left\langle \int_V \Phi((y+gA) \cap hB)\,dy, (g\kappa)  (h\lambda)
\right\rangle\, dg\, dh  \\
& = \int_G \int_{\overline{G}} \left\langle \Phi(\bar g A \cap hB), (g\kappa)  (h\lambda)
\right\rangle\, d\bar g\, dh \\ 
& = \int_{\overline{G}} \left\langle \Phi(A \cap \bar g^{-1} B), \kappa  g^{-1}\lambda
\right\rangle \,d\bar g. 
\end{align*}
In the last step, we used again the covariance of $\Phi$, the invariance properties of the Haar measure on $\overline{G}$ and 
the fact that $G$ is endowed with the Haar probability measure.

It follows that 
\begin{displaymath}
(\pd_m^{s_1} \otimes \pd_m^{s_2}) \circ k^G_{s_1,s_2}(\Phi)=m_G^* \circ \pd_m^{s_1+s_2}(\Phi),
\end{displaymath}
which finishes the proof. 
\endproof

\begin{Corollary} \label{cor_relation_kinformulas}
 Intersectional and additive kinematic formulas are related by the Alesker-Fourier transform in the following way:
\begin{displaymath}
 a^G=\left(\mathbb{F}^{-1} \otimes \mathbb{F}^{-1}\right) \circ k^G \circ \mathbb{F},
\end{displaymath}
or equivalently 
\begin{displaymath}
 k^G=\left(\mathbb{F} \otimes \mathbb{F}\right) \circ a^G \circ \mathbb{F}^{-1}.
\end{displaymath}
\end{Corollary}

\proof 
We claim that
\begin{displaymath}
\langle \pd_m^s \circ \mathbb{F} (\Phi_1),\mathbb{F} (\Phi_2)\rangle=(-1)^s \langle \pd_m^s (\Phi_1),\Phi_2\rangle,  
\end{displaymath}
whenever $\Phi_1,\Phi_2 \in \TVal^{s,G}$. 
Indeed, if $f,g \in \bigoplus_{i \equiv s \mod 2} \mathcal{H}^n_i$, then by Theorem \ref{mainthm_fourier} and Proposition \ref{prop_poincare_pairing} (whose proofs will be independent of this corollary),
\begin{displaymath}
 \mathbb{F}\mu_{k,f} \cdot \mathbb{F}\mu_{n-k,g}=(-1)^s \mu_{n-k, f} \cdot \mu_{k,g}=(-1)^s \mu_{k, f} \cdot \mu_{n-k,g}.
\end{displaymath}
Applying this formula to the components of $\Phi_1,\Phi_2$ and observing \eqref{eq_decomposition_sym} yields the claim.

Since Lemma \ref{lemma_prop_fourier} iv) implies $\mathbb{F}^2=(-1)^s\id$ on $\TVal^{s,G}$, we obtain that 
\begin{displaymath}
\langle \pd_m^s \circ \mathbb{F} (\Phi_1), \Phi_2\rangle= \langle \pd_m^s(\Phi_1),\mathbb{F} (\Phi_2)\rangle,  
\end{displaymath}
and hence
\begin{displaymath}
\mathbb{F}^* \circ \pd_m^s=\pd_m^s \circ \mathbb{F} .
\end{displaymath}
Using Lemma \ref{lemma_relation_pds}, we get 
\begin{displaymath}
 (\pd_m^s)^{-1} \circ \mathbb{F}^* \circ \pd_c^s = (-1)^s  (\pd_m^s)^{-1} \circ \mathbb{F}^* \circ \pd_m^s =  (-1)^s \mathbb{F}.
\end{displaymath}
Hence, the signs in the outer square in Theorem \ref{thm_ftaig} cancel out and the resulting commutative diagram is  
\begin{displaymath}
 \xymatrix{\TVal^{s_1+s_2,G} \ar[r]^<<<<<{a^G_{s_1,s_2}} \ar[d]_{\mathbb{F}} & \TVal^{s_1,G} \otimes
\TVal^{s_2,G}
\ar[d]_{\mathbb{F} \otimes \mathbb{F}} \\
\TVal^{s_1+s_2,G} \ar[r]^<<<<<{k^G_{s_1,s_2}}  & \TVal^{s_1,G}
\otimes \TVal^{s_2,G}
}
\end{displaymath}
which completes the proof. 
\endproof

In order to write down explicit formulas, we need to compute product and convolution of tensor valuations. This
will be achieved in the next two sections. 

\section{Alesker-Fourier transform of spherical valuations}
\label{sec_fourier}

The aim of this section is to prove Theorem \ref{mainthm_fourier}, i.e., to compute in an explicit way the Alesker-Fourier transform of spherical valuations. In Subsection \ref{subsec_fourier} we will restate this theorem in terms of tensor valuations. However, the connection between spherical valuations and 
tensor valuations will already be used to complete the proof of Theorem \ref{mainthm_fourier} in Section \ref{subsec_rel_tensor}.

\subsection{The Alesker-Fourier transform on translation invariant valuations}

The Fouri\-er transform on translation invariant valuations was introduced by Alesker (\cite{alesker03_un}, \cite{alesker_fourier}). Given an $n$-dimensional vector space $V$, the Alesker-Fourier transform is a map 
\begin{displaymath}
 \mathbb{F}: \Val^{sm}(V) \to \Val^{sm}(V^*) \otimes \Dens(V),
\end{displaymath}
where $\Dens(V)$ denotes the $1$-dimensional vector space of (complex) Lebesgue measures on $V$.

We will need the following properties of the Alesker-Fourier transform. 
\begin{Lemma} \label{lemma_prop_fourier}
\begin{enumerate}
\item $\mathbb{F}$ commutes with the action of $\GL(V)$. 
\item If $\phi \in \Val_k^{sm}(V)$, then $\mathbb{F}( \phi) \in \Val_{n-k}^{sm}(V^*) \otimes \Dens(V)$,
$k\in\{0,\ldots,n\}$.
\item Let $V$ be endowed with a Euclidean scalar product. Identifying $V^* \cong V$ and $\Dens(V) \cong \mathbb{C}$,
the Alesker-Fourier transform on {\it even} valuations satisfies 
\begin{displaymath}
 \Kl_{\mathbb{F}( \phi)}(E)=\Kl_\phi(E^\perp), \quad E \in \Gr_{n-k}(V), \phi \in \Val^{+,sm}_k(V),
\end{displaymath}
for $k\in\{0,\ldots,n\}$.
\item For $\phi \in \Val^{sm}$, we have the Plancherel-type formula
\begin{displaymath}
\mathbb{F}^2({\phi})(K)=\phi(-K),\quad K\in\mathcal{K}(V).
\end{displaymath}
\item For $\phi_1,\phi_2 \in \Val^{sm}$, we have 
\begin{displaymath}
 \mathbb{F}( \phi_1) *  \mathbb{F}( \phi_2)= \mathbb{F}(\phi_1 \cdot \phi_2).
\end{displaymath}
\end{enumerate}
\end{Lemma}

We will also use that the Alesker-Fourier transform is compatible with pull-backs, see Proposition \ref{prop_compatibility_pull_backs} for a precise statement.

\subsection{Representations of $\SO(n)$}

Let $V$ be a Euclidean vector space of dimension $n$ and put $\SO(n):=\SO(V)$. It is well-known that equivalence classes of complex irreducible (finite-dimensional) representations of $\SO(n)$ are
indexed by their highest weights. The possible highest weights are tuples $(\lambda_1,\lambda_2,\ldots,\lambda_{\lfloor
\frac{n}{2}\rfloor})$ of integers such that  
\begin{enumerate}
 \item 
$\lambda_1 \geq \lambda_2 \geq \ldots \geq \lambda_{\lfloor
\frac{n}{2}\rfloor} \geq 0$ if  $ n $ is odd,
\item $\lambda_1 \geq \lambda_2 \geq \ldots \geq |\lambda_{\frac{n}{2}}| \geq 0$ if $ n $ is  even.
\end{enumerate}

Given $\lambda=(\lambda_1,\ldots,\lambda_{\lfloor \frac{n}{2}\rfloor})$ satisfying this condition, we will denote
the corresponding equivalence class of representations by $\Gamma_\lambda$.

\begin{Proposition}
There is an isomorphism of $\SO(n)$-modules
\begin{displaymath}
 L^2(S^{n-1}) \cong \bigoplus_{s=0}^\infty \Gamma_{(s,0,\ldots,0)}.
\end{displaymath}
The $\SO(n)$-isotypical component in this decomposition contains precisely the spherical harmonics of degree $s$, i.e., 
$\Gamma_{(s,0,\ldots,0)} \cong \mathcal{H}^n_s$. Such representations are called {\it spherical}. 
\end{Proposition}

We refer to \cite{schneider_book14} for more information on spherical harmonics. In particular, we will need the following well-known fact. Let $\Sym^s V^*$ be the $\SO(n)$-module of all polynomials on $V$ of degree $s$. Then 
\begin{equation} \label{eq_decomposition_sym}
 \Sym^s V^* \cong \bigoplus_{i=0}^{\left\lfloor \frac{s}{2} \right\rfloor} \mathcal{H}_{s-2i}^n.
\end{equation}

\begin{Theorem}[\cite{strichartz75, takeuchi73}]
Let $\Gr_k$ denote the Grassmann manifold consisting of all $k$-dimensional subspaces in $\R^n$. The
$\SO(n)$-module $L^2(\Gr_k)$ decomposes as 
\begin{displaymath}
 L^2(\Gr_k) \cong \bigoplus_\lambda \Gamma_\lambda,
\end{displaymath}
where $\lambda$ ranges over all highest weights such that $\lambda_i=0$ for $i
> \min\{k,n-k\}$ and such that all $\lambda_i$ are even.  
\end{Theorem}

The decomposition of the $\SO(n)$-module $\Val_k$ has recently been obtained in
\cite{alesker_bernig_schuster}. 

\begin{Theorem}[\cite{alesker_bernig_schuster}]
There is an isomorphism of $\SO(n)$-modules
\begin{displaymath}
 \Val_k \cong \bigoplus_\lambda \Gamma_\lambda,
\end{displaymath}
where $\lambda$ ranges over all highest weights such that $|\lambda_2| \leq 2$, $|\lambda_i| \neq 1$ for all $i$ and
$\lambda_i=0$ for $i > \min\{k,n-k\}$. In particular, these decompositions are
multiplicity-free.
\end{Theorem}

\begin{Definition}
  Let $\Val_{k,s}$ for $ 1 \leq k \leq n-1$ and $ s \neq 1$ be
the $\SO(n)$-isotypical component of $\Val_k$ corresponding to $\Gamma_{(s,0,\ldots,0)}$. An element in $\bigoplus_{k,s} \Val_{k,s}$ is called {\it spherical valuation.} 
\end{Definition}

\begin{Corollary}\label{Corhns}
For $s\neq 1$, the map $\mathcal{H}^n_s \to \Val_{k,s}$, $f \mapsto \mu_{k,f}$, is an isomorphism of $\SO(n)$-modules. Moreover, 
\begin{equation} \label{eq_schur}
 \mathbb{F} (\mu_{k,f})=c_{n,k,s} \mu_{n-k,f}
\end{equation}
for some constant $c_{n,k,s} \in \mathbb{C}$ which only depends on $n,k,s$. 
\end{Corollary}

\proof
Since the non-trivial map $\mathcal{H}^n_s \to \Val_{k,s}, f \mapsto \mu_{k,f}$ commutes with $\SO(n)$, it is (by Schur's lemma) an isomorphism. 

The composition
\begin{displaymath}
 \mathcal{H}^n_s \cong \Val_{k,s} \stackrel{\mathbb{F}}{\longrightarrow} \Val_{n-k,s} \cong \mathcal{H}^n_s
\end{displaymath}
is an endomorphism of the irreducible $\SO(n)$-module $\mathcal{H}^n_s$ which commutes with $\SO(n)$. By Schur's lemma again,
it is given by multiplication by some scalar. 
\endproof

\subsection{Proof in the even case}
\label{subsec_proof_even}

We are going to prove the theorem in the even case,  i.e., evaluate the constants $c_{n,k,s}$ when $s$ is even. 
This case is simpler, and it will be useful to double check
the constants in the  more involved general case. 

For $k,l\in\{0,\ldots,n\}$ and $L\in \Gr_l$, we write $\Gr_k(L):=\{H\in \Gr_k:H\subset L\}$ if $k\le l$ and 
$\Gr_k(L):=\{H\in \Gr_k:L\subset H\}$ if $k\ge l$. 
 Then, for $0 \leq k,l \leq n$, we let $  R_{l,k}:C^\infty(\Gr_k) \to C^\infty(\Gr_l)$ denote the
Radon transform
\begin{displaymath}
 R_{l,k}f(L)=\int_{\Gr_k(L)} f(H)\, dH,
\end{displaymath}
where $dH$ is the Haar probability measure on $\Gr_k(L)$. Note that $R_{l,k}1=1$.

We will also need the map $\perp:C(\Gr_k) \to C(\Gr_{n-k})$ defined by $(\perp f)(U):=f(U^\perp)$ for 
$U\in \Gr_{n-k}$. 

With respect to the $L^2$-scalar products on $\Gr_k$ and $\Gr_l$ (all measures are normalized), we have 
\begin{displaymath}
 \langle R_{l,k}f,g\rangle=\langle f,R_{k,l}g\rangle,\quad  f \in C^\infty(\Gr_k), g \in C^\infty(\Gr_l).
\end{displaymath}

Furthermore, we have 
\begin{displaymath}
 R_{k,l}=\perp \circ R_{n-k,n-l} \circ \perp.
\end{displaymath}

Note also that 
\begin{displaymath}
 R_{m,k}=R_{m,l} \circ R_{l,k}
\end{displaymath}
whenever $k \leq l \leq m$ or $k \geq l \geq m$.

The next lemma is a simple case of a more general theorem by Grinberg \cite{grinberg86}. For the reader's convenience,
we give a direct proof. For an even function $f:S^{n-1}\to\R$, we define $f(U):=f(u)$ if $U$ is the linear span of $u\in S^{n-1}$. 
Hence $f$ can also be considered as a continuous function on $\Gr_1$. 

\begin{Lemma}\label{lemma_grin_specialcase}
Let $s\in\mathbb{N}_0$ be even, $f\in \mathcal{H}^n_s$ and $1\le k\le n-1$. Then
\begin{equation} \label{eq_iterated_radonspecial}
(R_{1,k} \circ R_{k,1})f=
\frac{\Gamma\left(\frac{n-1}{2}\right)\Gamma\left(\frac{k}{2}\right)
\Gamma\left(\frac{s+n-k}{2}\right)\Gamma\left(\frac{s+1}{2}\right)}{ \Gamma \left(
\frac12\right)
\Gamma\left(\frac{n-k}{2 }\right)\Gamma\left(\frac{s+k}{2}
\right)\Gamma\left(\frac{s+n-1}{2}\right)} f. 
\end{equation}
\end{Lemma}

\proof 
The case $k=1$ is trivial, so let us assume $2 \leq k \leq n-1$. For $v\in S^{n-1}$, let $[v] \in \Gr_1$ denote the span
of $v$. Then 
\begin{align*}
(R_{1,k} \circ R_{k,1})(f)([v])&=\int_{\Gr_k([v])}\int_{\Gr_1(L)}f(U)\, dU\, dL\\
&=\frac{1}{\omega_k}\int_{\Gr_k([v])}\int_{L\cap S^{n-1}}f(u)\, du\, dL,
\end{align*}
where $du$ denotes spherical Lebesgue measure on $L\cap S^{n-1}$. 

A special case of Theorem 7.2.2 in \cite{schneider_weil08} yields for an arbitrary continuous function $g$ on $\R^n$
that
\begin{displaymath}
\int_{B^n}g(x)\, dx=\frac{\omega_{n-1}}{\omega_{k-1}}\int_{\Gr_k([v])}\int_{L\cap B^n}g(y)|\det(y,v)|^{n-k}\, dy\,
dL.
\end{displaymath}
Here  $|\det(y,v)|=\|y\|\sqrt{1-\langle y/\|y\|,v\rangle^2}$ and $dy$ denotes Lebesgue measure on $L$.

Applying this to the homogeneous function $f$ and introducing polar coordinates, we get
\begin{displaymath}
(R_{1,k} \circ R_{k,1})(f)([v])=\frac{\omega_{k-1}}{\omega_{n-1}\omega_k} \int_{S^{n-1}} f(u)
\left(1-\langle u,v\rangle^2\right)^\frac{k-n}{2}\, du,
\end{displaymath}
where $du$ denotes spherical Lebesgue measure on $S^{n-1}$.
 
An application of the Funk-Hecke theorem to the integral on the right-hand side (see \cite[(1.10) and (1.21)]{rubin02}) yields
\begin{displaymath}
(R_{1,k} \circ R_{k,1})(f)([v])=\frac{\Gamma\left(\frac{n-1}{2}\right)\Gamma\left(\frac{k}{2}\right)
\Gamma\left(\frac{s+n-k}{2}\right)\Gamma\left(\frac{s+1}{2}\right)}{  \Gamma \left(
\frac12\right)
\Gamma\left(\frac{n-k}{2 }\right)\Gamma\left(\frac{s+k}{2}
\right)\Gamma\left(\frac{s+n-1}{2}\right)}
f(v).
\end{displaymath}
\endproof

We will also need the spherical Radon transform on even functions, which is given by 
\begin{displaymath}
 \perp \circ R_{n-1,1}: C_e^\infty(S^{n-1}) \to C_e^\infty(S^{n-1}).
\end{displaymath}
If $f$ is a harmonic polynomial of even degree $s$, then 
\begin{equation} \label{eq_spherical_radon}
( \perp \circ R_{n-1,1}) f=  (-1)^\frac{s}{2} \frac{\Gamma\left(\frac{n-1}{2}\right)
\Gamma\left(\frac{s+1}{2}\right)}{\Gamma\left(\frac12\right)
\Gamma\left(\frac{s+n-1}{2}\right)} f,
\end{equation}
compare \cite[Lemma 3.4.7]{groemer}.

\begin{Lemma} \label{lemma_klain_fct}
For $k\in \{1,\ldots,n-1\}$, the Klain function of $\mu_{k,f}$ is given by 
\begin{displaymath}
\Kl_{\mu_{k,f}}=(\perp \circ R_{n-k,1}) f . 
\end{displaymath}
\end{Lemma}

\proof
Let $E \in \Gr_k$, and let $B_E \subset E$ be the unit ball inside $E$. Then $S_k(B_E,\cdot)$ equals $\kappa_k
\omega_{n-k} \binom{n-1}{k}^{-1}$ times the normalized volume measure of the unit sphere in $E^\perp$.
Therefore 
\begin{displaymath}
\Kl_{\mu_{k,f}}(E)= \frac{1}{\vol_k B_E} \frac{1}{ \omega_{n-k}} \binom{n-1}{k} \int_{S^{n-1}} f\, dS_k(B_E,\cdot)=
(R_{n-k,1} f)(E^\perp). 
\end{displaymath}
\endproof

\proof[Proof of Theorem \ref{mainthm_fourier} in the even case]
Let $f \in \mathcal{H}^n_s$, where $s$ is even. By \eqref{eq_schur}, 
\begin{displaymath}
 \mathbb{F} (\mu_{k,f})=c_{n,k,s}\mu_{n-k,f}
\end{displaymath}
for some constant $c_{n,k,s}$. Then  Lemma \ref{lemma_prop_fourier} (iii) and Lemma \ref{lemma_klain_fct} imply that
\begin{displaymath}
  R_{n-k,1}f=c_{n,k,s} (\perp \circ R_{k,1})f.
\end{displaymath}
We apply $\perp \circ R_{n-1,n-k}=R_{1,k} \circ \perp$ to both sides and obtain
\begin{displaymath}
 (\perp \circ R_{n-1,1}) f = c_{n,k,s} (R_{1,k} \circ R_{k,1} )f.
\end{displaymath}
By Lemma \ref{lemma_grin_specialcase} and \eqref{eq_spherical_radon} we find 
\begin{displaymath}
 (-1)^\frac{s}{2} \frac{\Gamma\left(\frac{n-1}{2}\right)
\Gamma\left(\frac{s+1}{2}\right)}{ \Gamma \left(
\frac12\right) \Gamma\left(\frac{s+n-1}{2}\right)} f = c_{n,k,s}
\frac{\Gamma\left(\frac{n-1}{2}\right)\Gamma\left(\frac{k}{2}\right)
\Gamma\left(\frac{s+n-k}{2}\right)\Gamma\left(\frac{s+1}{2}\right)}{ \Gamma \left(
\frac12\right)
\Gamma\left(\frac{n-k}{2 }\right)\Gamma\left(\frac{s+k}{2}
\right)\Gamma\left(\frac{s+n-1}{2}\right)} f,
\end{displaymath}
and therefore 
\begin{displaymath}
 c_{n,k,s}= (-1)^\frac{s}{2}
\frac{\Gamma\left(\frac{n-k}{2}\right)\Gamma\left(\frac{s+k}{2}\right)}{\Gamma\left(\frac{k}{2}\right)\Gamma\left(\frac{
s+n-k}{2}\right)},
\end{displaymath}
which finishes the proof of Theorem \ref{mainthm_fourier} in the even case.  
\endproof

\subsection{The two-dimensional case}
\label{subsec_twodim}

The Alesker-Fourier transform on smooth valuations on a two-dimensional Euclidean vector space admits a simple description,
compare \cite{alesker_barcelona} which corrects a formula from \cite{alesker_fourier}. Using this description, we can
give a
direct proof of Theorem \ref{mainthm_fourier}, which will be useful as induction start in the general proof. 

Let $V \cong \R^2 \cong \mathbb{C}$ be a two-dimensional Euclidean vector space. Let 
\begin{displaymath}
 \mu_{1,f}(K):=\frac12 \int_{S^1} f(y) \, dS_1(K,y)
\end{displaymath}
with a smooth function $f:S^1 \to \mathbb{C}$. We may decompose $f$ as 
\begin{displaymath}
 f=f_++f_-^{\hol}+f_-^{\anti},
\end{displaymath}
where $f^+$ is the even part of $f$ and $f_-^{\hol}, f_-^{\anti}$ are the holomorphic and anti-holomorphic part of the
odd part of $f$, respectively. 

Then, according to \cite{alesker_barcelona}, the Alesker-Fourier transform satisfies 
\begin{displaymath}
 \mathbb{F}(\mu_{1,f})(K)=\frac12 \int_{S^1} (f_+(\i y)+f_-^{\hol}(\i y)-f_-^{\anti}(\i y))\,dS_1(K,y). 
\end{displaymath}

We claim that this formula is equivalent to the case $n=2, k=1$ of  Theorem \ref{mainthm_fourier}.  The
case $s=0$ is trivial. The space $\mathcal{H}_2^s$ for $s \geq 2$ is of dimension $2$, with the restrictions to the
circle of the complex valued
functions
$z^s, \bar z^s$ as basis. We consider four cases. 

\begin{itemize}
 \item Let $s$ be even and $f(z)=z^s$. Then $f_+(z)=z^s,f_-(z)=0$. Hence $f_+(\i y)=\i^s y^s=\i^s f_+(y)$ and
therefore
$\mathbb{F}(\mu_{1,f})=\i^s \mu_{1,f}$.
 \item Let $s$ be even and $f(z)=\bar z^s$. Then $f_+(z)=\bar z^s,f_-(z)=0$. Thus $f_+(\i y)=(-\i)^s \bar
y^s=\i^s
f_+(y)$ and therefore $\mathbb{F}(\mu_{1,f})=\i^s \mu_{1,f}$.
\item Let $s$ be odd and $f(z)=z^s$. Then $f_+(z)=0,f_-^{\hol}(z)=z^s, f_-^{\anti}(z)=0$. Hence $f_-^{\hol}(\i
y)=\i^s
y^s=\i^s f_-^{\hol}(y)$ and $\mathbb{F}(\mu_{1,f})=\i^s \mu_{1,f}$.
\item Let $s$ be odd and $f(z)=\bar z^s$. Then $f_+(z)=0,f_-^{\hol}(z)=0, f_-^{\anti}(z)=\bar z^s$. Hence 
$f_-^{\anti}(\i y)=(-\i)^s \bar y^s=- \i^s f_-^{\anti}(y)$ and $\mathbb{F}(\mu_{1,f})=\i^s \mu_{1,f}$.
\end{itemize}

In all cases, the result confirms the formula from Theorem \ref{mainthm_fourier}.

\subsection{Outline of the proof in the general case}
\label{subsec_outline}

For the remainder of the proof, we assume that $n \geq 3$. Our argument in this case is inspired by the computation of the spectrum of the $\cos^\lambda$-transform due to
\'Olafsson-Pasquale \cite{olafsson_pasquale}: we use the action of the larger group $\SL(n)$ to move between $K$-types
of
$\SO(n)$. 
Anticipating the results to be proved in Subsections \ref{subsec_laplacian}, \ref{subsec_poincare} and 
\ref{subsec_liealgebraaction}, we now deduce the recursion formula 
\begin{equation} \label{eq_quotient_fourier_s}
 \frac{c_{n,k,s+2}}{c_{n,k,s}}=- \frac{k+s}{n-k+s}.
\end{equation}
for the multipliers of the Alesker-Fourier transform. 

The first step in the proof of \eqref{eq_quotient_fourier_s} is Proposition \ref{prop_poincare_pairing}, where we will show  for $f,h \in C^\infty(S^{n-1})$ that 
\begin{equation}\label{firststep}
 \langle \pd_m(\mu_{k,f}),\mu_{n-k,h}\rangle=\langle \pd_m(\mu_{n-k,f}),\mu_{k,h}\rangle,
\end{equation}
which is equivalent to 
\begin{displaymath}
 \mu_{k,f} \cdot \mu_{n-k,h}=\mu_{n-k,f} \cdot \mu_{k,h}.
\end{displaymath}

Next, let $\mathfrak{g}$ be the Lie algebra of $G=\SL(n)$, i.e., the space of trace free square matrices. Then
$\mathfrak{g}=\mathfrak{so}(n) \oplus \mathfrak{h}$, where $\mathfrak{so}(n)$ is the Lie algebra of $\SO(n)$ (i.e., the
space of anti-symmetric  matrices) and $\mathfrak{h}$ is the space of trace free symmetric matrices. 

The actions of $G$ on $\Val(V)$ and $\Val(V^*)$ are given by $(g\mu)(K):=\mu(g^{-1}K)$ for $\mu \in \Val(V)$ and $K \in \mathcal{K}(V)$, and by  $(g\phi)(L):=\phi(g^*K)$ for $\phi \in \Val(V^*)$ and $L \in \mathcal{K}(V^*)$. These actions induce actions of the Lie algebra $\mathfrak{g}$ on the smooth subspaces $\Val^{sm}(V)$ and $\Val^{sm}(V^*)$ as follows. An element $X \in \mathfrak{g}$ may be represented as $g'(0)$ for some smooth curve $g:(-\epsilon,\epsilon) \to G$ with $g(0)=\mathrm{id}$. Then $X\mu:=\left.\frac{d}{dt}\right|_{t=0} g(t)\mu$ for $\mu \in \Val^{sm}(V)$, and similarly $X\phi:=\left.\frac{d}{dt}\right|_{t=0} g(t)\phi$ for $\phi \in \Val^{sm}(V^*)$.   

We will use the Euclidean scalar product to identify $V$ and $V^*$. The actions of $\mathfrak{h}$ on $\Val^{sm}(V)$ and on $\Val^{sm}(V^*) \cong \Val^{sm}(V)$ then differ by a sign. 

The Alesker-Fourier transform $\mathbb{F}:\Val^{sm}(V) \to \Val^{sm}(V^*) \cong \Val^{sm}(V)$ is continuous and commutes with the action of
$\mathfrak{h}$, and hence   
\begin{equation} \label{eq_commuting_f_x}
 \mathbb{F} (X \mu_{k,f})=-X \mathbb{F} (\mu_{k,f}),\quad X \in \mathfrak{h}.
\end{equation}

For $X\in \mathfrak{h}$, we define { $\xi_X \in C^\infty(S^{n-1})$ by}
\begin{displaymath}
 \xi_X(y):=\langle Xy,y\rangle, \quad y \in S^{n-1}. 
\end{displaymath}

We let $\pi_{k,s}$ denote the projection of $\Val_k$ onto its isotypical component of spherical type
$\Gamma_{(s,0,\ldots,0)}$. Since $\mathbb F$ commutes with $\SO(n)$, we have $\pi_{n-k,s} \circ \mathbb F=\mathbb F \circ \pi_{k,s}$. 
 
Let $h \in \mathcal{H}^n_{s+2}$. By Schur's lemma, 
\begin{displaymath}
 \langle \pd_m(\tau),\mu_{k,h}\rangle = \langle \pd_m \circ \pi_{n-k,s+2}(\tau),\mu_{k,h}\rangle
\end{displaymath}
for every $\tau \in \Val^{sm}_{n-k}$.

Given $f \in \mathcal{H}^n_s$, we may write $\pi_{k,s+2}(X \mu_{k,f})=\mu_{k,g}$ with a unique $g \in
\mathcal{H}^n_{s+2}$. Using \eqref{firststep}, we compute
\begin{align}
 \langle \pd_m \circ \mathbb{F}(X \mu_{k,f}),\mu_{k,h}\rangle & =  \langle \pd_m \circ \pi_{n-k,s+2} \circ \mathbb{F} (X
\mu_{k,f}),\mu_{k,h}\rangle\nonumber\\
 & = \langle \pd_m \circ \mathbb{F} \circ \pi_{k,s+2} (X
\mu_{k,f}),\mu_{k,h}\rangle\nonumber\\
& = \langle \pd_m \circ \mathbb{F} (\mu_{k,g}),\mu_{k,h}\rangle\nonumber\\
& = c_{n,k,s+2} \langle \pd_m(\mu_{n-k,g}),\mu_{k,h}\rangle\nonumber\\
& = c_{n,k,s+2} \langle \pd_m(\mu_{k,g}),\mu_{n-k,h}\rangle\nonumber\\
& = c_{n,k,s+2} \langle \pd_m(X\mu_{k,f}),\mu_{n-k,h}\rangle.\label{quote}
\end{align}On the other hand,
\begin{displaymath} 
 \langle \pd_m(X \mathbb{F} (\mu_{k,f})),\mu_{k,h}\rangle =  c_{n,k,s} \langle \pd_m(X \mu_{n-k,f}),\mu_{k,h}\rangle.
\end{displaymath}

By \eqref{eq_commuting_f_x} and \eqref{quote} it follows that  
\begin{equation}\label{quote2}
 c_{n,k,s+2} \langle \pd_m(X \mu_{k,f}),\mu_{n-k,h}\rangle=- c_{n,k,s} \langle \pd_m(X \mu_{n-k,f}),\mu_{k,h}\rangle.
\end{equation}

The core of the proof is the equation  
\begin{displaymath}
 \langle \pd_m(X\mu_{k,f}),\mu_{n-k,h}\rangle=\flag{n}{k} \frac{(-1)^{s+1}(s-1)(n+s+1)(n-k+s)}{\omega_n(n-1)}
\langle \xi_X \cdot f,h\rangle,
\end{displaymath}
which will be  established in Proposition \ref{zentraleProp}.
Here the scalar product on the right-hand side is the $L^2$-scalar product for functions on the sphere.  We
choose $X, f, h$ in such a way that $\langle \xi_X \cdot f,h\rangle \neq 0$ (we will see in Lemma \ref{nozero} that this is possible). 

Substituting this into \eqref{quote2}, we find the recursion equation \eqref{eq_quotient_fourier_s}. 

We may now use induction to prove \eqref{eq_fourier_constant}.
In the even case,  the induction start is $s=0$: since
$\mathbb{F}(\mu_k)=\mu_{n-k}$, we have $c_{n,k,0}=1$. 

In the odd case, the induction start will be $s=3$. In order to compute $c_{n,k,3}$, we use the Crofton formula from
\cite{hug_schneider_schuster_b} to compute the quotients $\frac{c_{n,k+1,3}}{c_{n,k,3}}$ (see Subsection 4.9.). This fixes all constants up to
a scaling which may depend on $n$. More precisely, 
\begin{equation} \label{eq_fourier_schur}
c_{n,k,s}=\epsilon_n \i^s \frac{\Gamma\left(\frac{n-k}{2}\right)\Gamma\left(\frac{s+k}{2}\right)}{\Gamma\left(\frac{k}{2}\right)\Gamma\left(\frac{
s+n-k}{2}\right)} \,,
\end{equation}
where $\epsilon_n$ depends only on $n$. Using functorial properties of the Alesker-Fourier transform and the proof in the
two-dimensional case, we  finally deduce 
that $\epsilon_n=1$ for all $n\ge 2$.   

\subsection{The Laplacian}
\label{subsec_laplacian}

The Laplacian acting on smooth functions on the sphere will be denoted by $\Delta_{S^{n-1}}$. This operator is self-adjoint with respect to the $L^2$-scalar product
\begin{displaymath}
 \langle f,g\rangle:=\int_{S^{n-1}} f g\, d\vol_{S^{n-1}}.
\end{displaymath}

We will also need the Laplacian on functions on $\R^n$ given by $\Delta_{\R^n}=-\sum_{i=1}^n \frac{\partial^2}{\partial
y_i^2}$. 

If $f$ is a smooth function which is defined in a neighborhood of the unit sphere, then, on the sphere $S^{n-1}$,
\begin{equation} \label{eq_relation_laplacians}
 \Delta_{S^{n-1}}f=\Delta_{\R^n}f+(n-1)\frac{\partial}{\partial r} f +\frac{\partial^2}{\partial r^2}f.
\end{equation}
Here $\frac{\partial}{\partial r}$ is the radial derivation.  We refer to (\cite{gallot_hulin_lafontaine}, 4.E) 
for more information. In particular, if $f$ is homogeneous of degree $s$, then \eqref{eq_relation_laplacians} 
implies that 
\begin{equation}\label{relmod}
\Delta_{S^{n-1}}f=\Delta_{\R^n}f+s(n+s-2)f \qquad\text{on }S^{n-1}.
\end{equation}
A homogeneous polynomial $p$ of degree $s$ with  $\Delta_{\R^n} p=0$ is called  {\it harmonic}. From
the previous equation, it follows that a harmonic polynomial of degree $s$ is an eigenfunction of the spherical
Laplacian with eigenvalue $s(n+s-2)$.

The natural operation of $\mathfrak{h}$ on $C^\infty(\R^n)$ is given by 
\begin{displaymath}
 \hat Xh:= \sum_{i,j} X_{ij}y_j \frac{\partial}{\partial y_i} h \in C^\infty(\R^n),
\end{displaymath}
where $h \in C^\infty(\R^n)$ and $X=(X_{ij})_{i,j=1}^n \in \mathfrak{h}$. 
We will write $h_i:=\frac{\partial}{\partial y_i} h$ and $h_{ij}:=\frac{\partial^2}{\partial y_i \partial y_j} h$.
It will be crucial that $\hat X$ maps homogeneous polynomials of degree $s$ to homogeneous polynomials of degree $s$.

\begin{Lemma} \label{lemma_laplacian_hatx}
Let $h$ be a homogeneous harmonic polynomial and $X \in \mathfrak{h}$. Then 
\begin{equation} 
[\Delta_{S^{n-1}}, \hat X](h)=-2\sum_{i,j=1}^n X_{ij} h_{ij}\qquad \text{on }S^{n-1}.
\end{equation}
\end{Lemma}

\proof 
Since $h$ is harmonic, we get
\begin{equation}\label{eqoben}
 \Delta_{\R^n}(\hat Xh)  = - \sum_{i,j=1}^n\sum_{l=1}^n X_{ij} \frac{\partial^2}{\partial y_l^2} \left(y_j h_i
\right) = -2 \sum_{i,j=1}^n X_{ij} h_{ij}.
\end{equation}
If $h$ is homogeneous of degree $s$, then so is $\hat Xh$. Hence \eqref{relmod} implies that
\begin{align*}
[\Delta_{S^{n-1}}, \hat X](h)&=\Delta_{S^{n-1}}(\hat Xh)-\hat X(\Delta_{S^{n-1}} h)\\
&=\Delta_{\R^n}(\hat X h)+s(n-2+s)\hat X h -\hat X(\Delta_{\R^n} h+s(n-2+s)h)\\
&=\Delta_{\R^n}(\hat X h),
\end{align*}
and the assertion follows from \eqref{eqoben}.
\endproof

\subsection{Poincar\'e pairing}
\label{subsec_poincare}

Let $V$ be an oriented $n$-dimensional Euclidean vector space. We use Euclidean coordinates $(x_1,x_2,\ldots,x_n)$. On a
second copy of $V$, we use coordinates $(y_1,\ldots,y_n)$. Differential forms on the sphere bundle $SV=\R^n \times
S^{n-1}$ will be written as restrictions of forms on $\R^n \times \R^n$. The exterior differential will be denoted by
$d$; the contraction of a form $\omega$ with a vector field $X$ will be denoted by $i_X \omega$. For an ordered
set $I=(i_1,\ldots,i_k)$ 
of indices, we use the abbreviation $dx_I:=dx_{i_1} \wedge \cdots \wedge dx_{i_k}$, and
similarly for $dy_I$. 

Let $\rho_k \in \Omega^{k,n-k-1}(SV)$, for $k\in\{0,\ldots,n-1\}$, be the Lipschitz-Killing form defined by 
\begin{displaymath}
\rho_k = \frac{1}{k!(n-k-1)!\omega_{n-k}} \sum_{\pi \in \mathcal{S}_n} \sgn(\pi) y_{\pi(1)}dy_{\pi(2) \ldots
\pi(n-k)}\wedge
dx_{\pi(n-k+1) \ldots \pi(n)}.
\end{displaymath}
Here $\pi$ ranges over the set $\mathcal{S}_n$ of all permutations of the elements $\{1,2,\ldots,n\}$ and $\sgn(\pi)$ is
the sign of the permutation $\pi$. Moreover, for $k=0$ we interpret $dx_{\pi(n-k+1) \ldots \pi(n)}$ as $1$ and if
$k=n-1$ then $dy_{\pi(2) \ldots
\pi(n-k)}$ is taken as $1$.

This definition is independent of the choice of a positively oriented orthonormal basis of $\R^n$. Hence, in the
pointwise
computations which will
follow, we will usually adapt our coordinate system in such a way that we only have to consider the point
$y= e_1=(1,0,\ldots,0)$. 
At that point we have $dy_1=0$. Moreover, the contact form $\alpha$ and  its differential $d\alpha$ are given by  
\begin{displaymath}
\alpha=dx_1,\quad d\alpha =\sum_{k=2}^n dy_k \wedge dx_k.
\end{displaymath}

Recall that $\nc(K)$ denotes the normal cycle (normal bundle) of a compact convex set $K$ in $V$. Hence,
for $h \in C^\infty(S^{n-1})$ 
and $k\in\{0,\ldots,n-1\}$, we have
\begin{displaymath}
\mu_{k,h}(K)=\int_{\nc(K)} h \rho_k.
\end{displaymath}

Let us briefly recall the Rumin differential $D:\Omega^{n-1}(S\R^n) \to \Omega^n(S\R^n)$. A multiple of the contact form
$\alpha$ is called {\it vertical}. For instance, it is easily checked by a direct calculation that $d\rho_k$ is vertical for all $k$ (see also Lemma 3.1 in \cite{Fu98}).
Given an $(n-1)$-form $\gamma$, there exists a unique vertical form, which we
write
as $\alpha \wedge \xi$, such that $d(\gamma+\alpha \wedge \xi)$ is again vertical. The operator
$D\gamma:=d(\gamma+\alpha \wedge \xi)$ is called {\it Rumin operator}. A similar definition makes sense on any
contact
manifold. We refer to \cite{rumin94} for the basic properties of $D$ and to \cite{alesker_bernig, bernig_quat09,
bernig_fu06} for applications in integral geometry, some of which will be explained in more detail below.  

\begin{Proposition}  \label{prop_rumin}
 Let $h \in C^\infty(S^{n-1})$. For $k\in\{0,\ldots,n-2\}$, define an $(n-2)$-form $\xi_{k,h}$ by 
\begin{displaymath}
\xi_{k,h}:= i_{\grad h} \rho_k \in \Omega^{k,n-k-2}(S\R^n).
\end{displaymath}
Then, for $k\in\{1,\ldots,n-1\}$,
\begin{displaymath}
 D(h \rho_k)=d\left(h\rho_k-\frac{\omega_{n-k+1}}{\omega_{n-k}} \alpha \wedge \xi_{k-1,h}\right)=h
d\rho_k+\frac{\omega_{n-k+1}}{\omega_{n-k}} \alpha \wedge d\xi_{k-1,h}.
\end{displaymath}
At the point $(0,e_1)$, we have 
\begin{align}
 D(h \rho_k) & =  \frac{1}{(k-1)!(n-k-1)!\omega_{n-k}} \alpha \wedge \sum_{\pi \in \mathcal{S}'_n} \sgn(\pi)
\left(h-h_1+h_{\pi(2)\pi(2)} \right)\cdot\nonumber\\
& \qquad\qquad \cdot dy_{\pi(2) \ldots
\pi(n-k+1)}\wedge
dx_{\pi(n-k+2) \ldots \pi(n)}\nonumber \\
& \quad +  \frac{(-1)^{n-k+1}(k-1)}{(k-1)!(n-k-1)!\omega_{n-k}} 
\alpha \wedge \sum_{\pi \in \mathcal{S}_n'} \sgn(\pi) h_{\pi(2)\pi(n-k+2)} \cdot\nonumber \\
& \qquad\qquad \cdot dy_{\pi(3)\ldots
\pi(n-k+2)} \wedge dx_{\pi(n-k+2) \ldots \pi(n)}, \label{eq_rumin_hrho}
\end{align}
where $\mathcal{S}_n'$ is the set of permutations $\pi \in \mathcal{S}_n$ such that $\pi(1)=1$. 
\end{Proposition}

\proof
Let $k\in\{1,\ldots,n-1\}$. We claim that 
\begin{displaymath}
 \omega_{n-k}dh \wedge \rho_k=\omega_{n-k+1} d\alpha \wedge \xi_{k-1,h},
\end{displaymath}
from which the first part of the proposition follows.

By $\On(n)$-equivariance, it is enough to do the computation at the point $y= e_1$. We compute 
\begin{align*}
 dh \wedge \rho_k & = \frac{1}{k!(n-k-1)!\omega_{n-k}} \sum_{i=2}^n h_i dy_i \wedge \sum_{\pi \in \mathcal{S}_n'}
\sgn(\pi)\\
& \quad 
dy_{\pi(2) \ldots \pi(n-k)} \wedge
dx_{\pi(n-k+1) \ldots \pi(n)} \\
& = \frac{1}{k!(n-k-1)!\omega_{n-k}} \sum_{\pi \in \mathcal{S}_n'} \sgn(\pi)\\
& \quad \sum_{i=n-k+1}^n  h_{\pi(i)} dy_{\pi(i) \pi(2) \ldots \pi(n-k)}\wedge
dx_{\pi(n-k+1) \ldots \pi(n)} \\
& =  \frac{1}{(k-1)!(n-k-1)!\omega_{n-k}} \sum_{\pi \in \mathcal{S}_n'} \sgn(\pi)\\
& \quad 
h_{\pi(n-k+1)} dy_{\pi(n-k+1) \pi(2) \ldots \pi(n-k)}\wedge
dx_{\pi(n-k+1) \ldots \pi(n)} \\
& =   \frac{(-1)^{n-k-1}}{(k-1)!(n-k-1)!\omega_{n-k}} \sum_{\pi \in \mathcal{S}_n'} \sgn(\pi)\\
& \quad 
h_{\pi(n-k+1)} dy_{\pi(2) \ldots \pi(n-k+1)} \wedge
dx_{\pi(n-k+1) \ldots \pi(n)}.
\end{align*}

Since 
\begin{align*}
\xi_{k-1,h} & =\frac{1}{(k-1)!(n-k)!\omega_{n-k+1}}\\
& \quad i_{\sum_{j=2}^n h_j
\frac{\partial}{\partial y_j}} \left(\sum_{\pi \in \mathcal{S}_n}
\sgn(\pi) y_{\pi(1)} dy_{\pi(2) \ldots \pi(n-k+1)}\wedge
dx_{\pi(n-k+2) \ldots \pi(n)}\right)\\
& =\frac{1}{(k-1)!(n-k-1)!\omega_{n-k+1}} \\
& \quad \sum_{\pi \in \mathcal{S}_n} 
\sgn(\pi) y_{\pi(1)} h_{\pi(2)} dy_{\pi(3) \ldots \pi(n-k+1)}\wedge
dx_{\pi(n-k+2) \ldots \pi(n)}
\end{align*}
we obtain that  
\begin{align*}
 d\alpha \wedge \xi_{k-1,h} & =\frac{1}{(k-1)!(n-k-1)!\omega_{n-k+1}}\sum_{\pi \in \mathcal{S}_n'}
\sgn(\pi) \\
& \quad \sum_{i=2}^n dy_{\pi(i)} \wedge dx_{\pi(i)} \wedge h_{\pi(2)} dy_{\pi(3) \ldots
\pi(n-k+1)}\wedge
dx_{\pi(n-k+2) \ldots \pi(n)}\\
& =\frac{1}{(k-1)!(n-k-1)!\omega_{n-k+1}} \sum_{\pi \in \mathcal{S}_n'}
\sgn(\pi) \\
& \quad dy_{\pi(2)} \wedge dx_{\pi(2)} \wedge h_{\pi(2)} dy_{\pi(3) \ldots
\pi(n-k+1)}\wedge
dx_{\pi(n-k+2) \ldots \pi(n)}\\
& =\frac{(-1)^{n-k-1}}{(k-1)!(n-k-1)!\omega_{n-k+1}}\sum_{\pi \in \mathcal{S}_n'}
\sgn(\pi) \\
& \quad h_{\pi(n-k+1)} dy_{\pi(2) \ldots \pi(n-k+1)} \wedge dx_{\pi(n-k+1) \ldots \pi(n)}.
\end{align*}

This proves the claim and thus the first part of the proposition.  The second part follows by an explicit
calculation of $d\rho_k$ and $d\xi_{k-1,h}$. 
\endproof

\begin{Proposition} \label{prop_poincare_pairing}
 Let $h \in \mathcal{H}^n_{s_1}, g \in \mathcal{H}^n_{s_2}$ with $s_1,s_2 \neq 1$. Then, for $1 \leq k \leq n-1$, 
\begin{displaymath}
 \mu_{k,h} \cdot \mu_{n-k,g}=\begin{cases} \frac{(-1)^{s_1}}{\omega_n}\flag{n}{k} \left\langle
h-\frac{1}{n-1}\Delta_{S^{n-1}} h,g\right\rangle \vol_n, & s_1=s_2, \\
0, & s_1 \neq s_2.
\end{cases}
\end{displaymath}
 
In particular, for $h,g\in C^\infty(S^{n-1})$,
\begin{displaymath}
 \langle \pd_m(\mu_{k,h}),\mu_{n-k,g}\rangle=\langle \pd_m(\mu_{n-k,h}),\mu_{k,g}\rangle.
\end{displaymath}
\end{Proposition}

\proof
Let us recall a formula from \cite{bernig_quat09}. Let $\pi:S\R^n \to \R^n$ be the projection and $\pi_*:\Omega^*(S\R^n)
\to \Omega^{*-n+1}(\R^n)$ be the corresponding fiber integration. If $\rho$ is a translation invariant $(2n-1)$-form on
$S\R^n$, then $\pi_*(\rho)$ is a multiple of the Lebesgue measure, i.e., $\pi_*\rho=c \vol_{\R^n}$ for some constant
$c$. We set $\int \rho:=c$. 

Given translation invariant valuations $\phi_1,\phi_2$, represented by differential forms $\gamma_1,\gamma_2 \in
\Omega^{n-1}(S\R^n)$, we have, by \cite[Theorem 1.3]{bernig_quat09}, 
\begin{displaymath}
\langle \pd_m(\phi_1),\sigma \phi_2\rangle=(-1)^n \int \gamma_1 \wedge D\gamma_2.
\end{displaymath}

We apply this formula with $\phi_1:=\mu_{n-k,g}, \phi_2:=\mu_{k,h}$. By the definition of the Verdier involution, we have $\sigma \mu_{k,h}=(-1)^{k+s_1}\mu_{k,h}$. 

The statement of the proposition thus follows if we can show that
\begin{equation} \label{eq_rho_Dh}
  \rho_{n-k} \wedge D(h \rho_{k})=\flag{n}{k} \frac{(-1)^{n+k}}{\omega_n} \left(h
-\frac{1}{n-1}\Delta_{S^{n-1}} h\right) \vol_{S\R^n}. 
\end{equation}

From \eqref{eq_rumin_hrho} we get 
\begin{align*}
 \rho_{n-k} \wedge D(h\rho_k) & = \frac{1}{(k-1)!(n-k)!\omega_k} \sum_{\pi \in \mathcal{S}'_n}
\sgn(\pi)
dy_{\pi(2) \ldots
\pi(k)}\wedge
dx_{\pi(k+1) \ldots \pi(n)} \\
& \quad \wedge \alpha \wedge \frac{1}{(k-1)!(n-k-1)!\omega_{n-k}} \sum_{\tau \in \mathcal{S}'_n} \sgn(\tau)\\
& \quad 
\left(h-h_1+h_{\tau(2)\tau(2)} \right) dy_{\tau(2) \ldots
\tau(n-k+1)} \wedge
dx_{\tau(n-k+2) \ldots \tau(n)}.
\end{align*}
For a fixed $\tau$, the wedge product will vanish except when $\{\pi(2),\ldots,\pi(k)\}$
and $\{\tau(2),\ldots,\tau(n-k+1)\}$ are disjoint. There are $(k-1)!(n-k)!$ such permutations. 

We thus get 
\begin{align*}
 &\rho_{n-k} \wedge D(h\rho_k)\\ 
 & = (-1)^{n+k} \frac{1}{(k-1)!(n-k-1)!\omega_k \omega_{n-k}} \sum_{\tau \in
\mathcal{S}'_n} 
\left(h-h_1+h_{\tau(2)\tau(2)} \right) \vol_{S\R^n}\\
& = (-1)^{n+k} \frac{1}{\omega_n}\flag{n}{k}
\left(h-\frac{\partial h}{\partial r}+\frac{1}{n-1} h_{22}+\ldots+\frac{1}{n-1}
h_{nn} \right) \vol_{S\R^n}\\
& = \frac{(-1)^{n+k}}{\omega_n}\flag{n}{k} h\left(1-  s_1-\frac{s_1(s_1-1)}{n-1}\right) \vol_{S\R^n} \\
& = \frac{(-1)^{n+k}}{\omega_n}\flag{n}{k} \left(h-\frac{\Delta_{S^{n-1}} h}{n-1}\right) \vol_{S\R^n},
\end{align*}
which proves \eqref{eq_rho_Dh}.
\endproof

\subsection{The Lie algebra action}
\label{subsec_liealgebraaction}

The  action of $G=\SL(n)$ on $\R^n$ induces an action by contactomorphisms on the sphere bundle $S\R^n=\R^n
\times S^{n-1}$ which is given by
\begin{displaymath}
\tau(g)(x,u)=\left(gx, \frac{g^{-\top}u}{\|g^{-\top}u\|}\right). 
\end{displaymath}
This action is compatible with the map
\begin{displaymath}
\Omega^{n-1}(S\R^n) \to \Val(\R^n),\quad \gamma \mapsto [K \mapsto
\int_{\nc(K)} \gamma].
\end{displaymath}
The corresponding action of the  subspace $\mathfrak{h} \subset \mathfrak{g}$ on $\Omega^*(S\R^n)$ will be
denoted by $(X,\gamma) \mapsto X \gamma$.  Explicitly, if $g:(-\epsilon,\epsilon) \to G$ is a smooth curve with $g(0)=\mathrm{id}$ and $g'(0)=X$, then $X \gamma=\left.\frac{d}{dt}\right|_{t=0} (\tau (g(t))^* \gamma$ for $\gamma \in \Omega^*(S\R^n)$.

In particular, $X$ acts  on $C^\infty(S\R^n)=\Omega^0(S\R^n)$ as
\begin{equation}\label{actiona}
 X x_i=\sum_j X_{ij}x_j, \quad X y_i=-\sum_j X_{ji}y_j+\underbrace{\left(\sum_{kl} X_{lk}y_ly_k\right)}_{=\xi_X(y)}
y_i. 
\end{equation}

Note that if $f$ is the restriction of a function $\hat f$ on $\R^n$, then
\begin{equation} \label{eq_relation_x_hatx}
 X f=-\hat X \hat f \vert_{S^{n-1}}+\xi_X \left. \frac{\partial}{\partial r}\right\vert_{r=1} \hat f.
\end{equation}

If $f,h \in C^\infty(S^{n-1})$, then 
\begin{equation} \label{eq_x_on_product}
 0=\langle Xf,h\rangle+\langle f,Xh\rangle+n \langle f,\xi_X h\rangle. 
\end{equation}
Indeed, let $\omega:=f h \vol_{S^{n-1}} \in \Omega^{n-1}(S^{n-1})$. Since the action of $\SL(n)$ on $S^{n-1}$ is
by orientation preserving diffeomorphisms, we have $\int_{S^{n-1}} g^*\omega=\int_{S^{n-1}} \omega$. Taking the
derivative yields 
\begin{displaymath}
 \int_{S^{n-1}} X\omega=0, \quad X \in \mathfrak{h}. 
\end{displaymath}
Since 
\begin{displaymath}
X\omega=X(f) h \vol_{S^{n-1}}+f X(h) \vol_{S^{n-1}}+ fh n \xi_X \vol_{S^{n-1}},
\end{displaymath}
equation \eqref{eq_x_on_product} follows. 
 
\begin{Proposition} \label{prop_action_liealgebra}
 Let $X \in \mathfrak{h}$ and $1 \leq k \leq n-1$. Then, at the special point $y= e_1$, 
\begin{align}
X\rho_k & \equiv \frac{ (n-k-2)}{k!(n-k-1)!\omega_{n-k}} \xi_X \sum_{\pi \in \mathcal{S}_n'} \sgn(\pi)
dy_{\pi(2) \ldots\pi(n-k)} \wedge dx_{\pi(n-k+1)\ldots \pi(n)} \label{X_rho1} \\
& \quad  -\frac{2(n-k-1)}{k!(n-k-1)!\omega_{n-k}} \sum_{\pi \in \mathcal{S}_n'} \sgn(\pi)
X_{\pi(2)\pi(2)}\nonumber\\
& 
\qquad\qquad dy_{\pi(2) \ldots\pi(n-k)} \wedge dx_{\pi(n-k+1)\ldots \pi(n)} \label{X_rho2}\\
& \quad + \frac{2(-1)^{n-k+1}(n-k-1)}{(k-1)!(n-k-1)!\omega_{n-k}} \sum_{\pi \in \mathcal{S}_n'} \sgn(\pi) X_{\pi(2)
\pi(n-k+1)} \nonumber \\
& \qquad\qquad  
dy_{\pi(3) \ldots \pi(n-k+1)} \wedge dx_{\pi(n-k+1)\ldots
\pi(n)} \mod \alpha. \label{X_rho3}
\end{align}
\end{Proposition}

\proof 
The form $\rho_k$ has three different types of factors: $y_{\pi(1)}$, the $(n-k-1)$ factors $dy_i$, and the $k$ factors $dx_j$. 
Correspondingly, $X\rho_k$ splits  as a sum of three terms
$A_1,A_2,A_3$. Using \eqref{actiona}, we get 
\begin{align*}
 A_1 & = - \frac{1}{k!(n-k-1)!\omega_{n-k}} \sum_{\pi,j} \sgn(\pi) X_{j\pi(1)}y_j
dy_{\pi(2)\ldots\pi(n-k)} \wedge dx_{\pi(n-k+1)\ldots \pi(n)}\\
& \quad + \frac{1}{k!(n-k-1)!\omega_{n-k}} \xi_X(y) \sum_{\pi} \sgn(\pi) y_{\pi(1)}
dy_{\pi(2)\ldots\pi(n-k)} \wedge dx_{\pi(n-k+1)\ldots \pi(n)},
\end{align*}
where the summation extends over $\pi\in\mathcal{S}_n$ and $j\in\{1,\ldots,n\}$. 
We evaluate this at $y=e_1$. Then $X_{11}=\xi_X(y)$. In the first sum, all terms with $j
\neq 1$ vanish, while
all terms with $\pi(1) \neq 1$ are divisible by $dy_1=0$ or by $dx_1=\alpha$. The terms with $j=1, \pi(1)=1$ cancel with
the second sum. Hence, at the special point $y=e_1$ we have 
\begin{displaymath}
A_1 \equiv 0 \mod \alpha.
\end{displaymath}

Next, 
\begin{align*}
 A_2 & = -\frac{n-k-1}{k!(n-k-1)!\omega_{n-k}} \sum_{\pi,j} \sgn(\pi) y_{\pi(1)} X_{j\pi(2)} dy_{j\pi(3)
\ldots\pi(n-k)} \wedge dx_{\pi(n-k+1)\ldots
\pi(n)}\\
& \quad +\frac{n-k-1}{k!(n-k-1)!\omega_{n-k}} \sum_{\pi} \sgn(\pi) y_{\pi(1)}\\
&\qquad\qquad d [\xi_X(y) y_{\pi(2)}]  dy_{\pi(3)
\ldots\pi(n-k)} \wedge dx_{\pi(n-k+1)\ldots
\pi(n)}\\
& = -\frac{n-k-1}{k!(n-k-1)!\omega_{n-k}} \sum_{\pi \in \mathcal{S}_n'} \sgn(\pi)
X_{\pi(2)\pi(2)}\\
&\qquad\qquad 
dy_{\pi(2)\pi(3) \ldots\pi(n-k)} \wedge dx_{\pi(n-k+1)\ldots \pi(n)}\\
& \quad +\frac{(-1)^{n-k+1}(n-k-1)}{(k-1)!(n-k-1)!\omega_{n-k}}  \sum_{\pi \in \mathcal{S}_n'} \sgn(\pi)
X_{\pi(n-k+1)\pi(2)} \\
&\qquad\qquad dy_{\pi(3)
\ldots\pi(n-k+1)} \wedge dx_{\pi(n-k+1)\ldots
\pi(n)}\\
& \quad +\frac{(n-k-1) \xi_X(y)}{k!(n-k-1)!\omega_{n-k}}  \sum_{\pi \in \mathcal{S}_n'} \sgn(\pi) dy_{\pi(2)
\ldots\pi(n-k)} \wedge dx_{\pi(n-k+1)\ldots
\pi(n)}.
\end{align*}

Finally, 
\begin{align*}
 A_3 & = \frac{1}{(k-1)!(n-k-1)!\omega_{n-k}} \sum_{\pi,j} \sgn(\pi) y_{\pi(1)} X_{\pi(n-k+1)j} \\
 &\qquad\qquad dy_{\pi(2)
\ldots\pi(n-k)} \wedge dx_{j \pi(n-k+2)\ldots
\pi(n)}\\
& \equiv \frac{1}{(k-1)!(n-k-1)!\omega_{n-k}} \sum_{\pi \in \mathcal{S}_n'} \sgn(\pi)
X_{\pi(n-k+1)\pi(n-k+1)} \\
 &\qquad\qquad dy_{\pi(2) \ldots\pi(n-k)} \wedge dx_{\pi(n-k+1) \ldots
\pi(n)}\\
& \quad + \frac{(-1)^{n-k+1}(n-k-1)}{(k-1)!(n-k-1)!\omega_{n-k}}  \sum_{\pi \in \mathcal{S}_n'} \sgn(\pi) X_{\pi(2)
\pi(n-k+1)}\\
 &\qquad\qquad 
dy_{\pi(3) \ldots \pi(n-k+1)} \wedge dx_{\pi(n-k+1)\ldots
\pi(n)} \mod\alpha.
\end{align*}

Since the trace of $X$ is zero, it follows that 
\begin{align*}
 X\rho_k & = A_1+A_2+A_3 \\
& \equiv \frac{\xi_X(y)(n-k-2)}{k!(n-k-1)!\omega_{n-k}}  \sum_{\pi \in \mathcal{S}_n'}
\sgn(\pi)
dy_{\pi(2) \ldots\pi(n-k)} \wedge dx_{\pi(n-k+1)\ldots \pi(n)} \\
& \quad  -\frac{2(n-k-1)}{k!(n-k-1)!\omega_{n-k}} \sum_{\pi \in \mathcal{S}_n'} \sgn(\pi)
X_{\pi(2)\pi(2)}\\
 &\qquad\qquad
dy_{\pi(2)\pi(3) \ldots\pi(n-k)} \wedge dx_{\pi(n-k+1)\ldots \pi(n)} \\
& \quad + \frac{2(-1)^{n-k+1}(n-k-1)}{(k-1)!(n-k-1)!\omega_{n-k}}  \sum_{\pi \in
\mathcal{S}_n'} \sgn(\pi) 
X_{\pi(2)
\pi(n-k+1)}\\
 &\qquad\qquad
dy_{\pi(3) \ldots \pi(n-k+1)} \wedge dx_{\pi(n-k+1)\ldots
\pi(n)} \mod \alpha,
\end{align*}
which yields the assertion of the proposition.
\endproof

\begin{Proposition} \label{prop_action_poincare}
Let $h \in \mathcal{H}^n_{s+2}$ and $1 \leq k \leq n-1$. Then 
 \begin{displaymath}
  X\rho_k \wedge D(h\rho_{n-k})=\left(a\, \Delta_{S^{n-1}}(\hat X h)+b\, \hat X h + c \,\xi_X { \cdot} h\right)
\Omega,
 \end{displaymath}
where 
\begin{align*}
 \Omega & :=\flag{n}{k} \frac{1}{\omega_n} \vol_{S\R^n},\\
 a & := \frac{(-1)^{k+1}(n-k-1)}{(n-1)(n-2)}, \\
 b & := \frac{(-1)^k(n-k-1)((n+s)(s+2)-4(s+1))}{(n-1)(n-2)},  \\
 c & := (-1)^{k+1}\frac{s+1}{n-1}\left((n-k-2)(n+s+3)+2\right).
\end{align*}
\end{Proposition}

\proof
Let $h$ be a harmonic polynomial of degree $s+2$. According to Proposition \ref{prop_rumin} we have 
\begin{align}
 &D(h \rho_{n-k})\nonumber \\
 & =  - \frac{(1+s)h}{(k-1)!(n-k-1)!\omega_k} \alpha \wedge \sum_{\pi \in \mathcal{S}'_n} \sgn(\pi)
dy_{\pi(2) \ldots
\pi(k+1)}\wedge
dx_{\pi(k+2) \ldots \pi(n)} \label{D_h_rho1}\\
& \qquad +  \frac{1}{(k-1)!(n-k-1)!\omega_k} \alpha \wedge \sum_{\pi \in \mathcal{S}'_n} \sgn(\pi)
h_{\pi(2)\pi(2)} \nonumber\\&\qquad\qquad 
dy_{\pi(2) \ldots
\pi(k+1)}\wedge
dx_{\pi(k+2) \ldots \pi(n)}  \label{D_h_rho2}\\
& \qquad +  \frac{(-1)^{k+1}(n-k-1)}{(k-1)!(n-k-1)!\omega_k} 
 \nonumber \\
& \qquad \qquad\alpha \wedge \sum_{\pi \in \mathcal{S}_n'} \sgn(\pi) h_{\pi(2)\pi(k+2)} dy_{\pi(3)\ldots
\pi(k+2)} \wedge dx_{\pi(k+2) \ldots \pi(n)}. \label{D_h_rho3}
\end{align}

The non-zero terms in $X\rho_k \wedge D(h\rho_{n-k})$ are given by 
\begin{align*}
 \eqref{X_rho1} \wedge \eqref{D_h_rho1} & = (-1)^{k+1} (n-k-2)(s+1)  h \xi_X \Omega,\\
\eqref{X_rho1} \wedge \eqref{D_h_rho2} & = (-1)^{k+1}  \frac{(n-k-2)(s+2)(s+1)}{n-1} h \xi_X \Omega,\\
 \eqref{X_rho2} \wedge \eqref{D_h_rho1} & = (-1)^{k+1} \frac{2(n-k-1) (s+1)}{n-1} h \xi_X \Omega, \\
\eqref{X_rho2} \wedge \eqref{D_h_rho2} & = \frac{2(-1)^k(n-k-1)}{(n-1)(n-2)} \sum_{j=2}^n X_{jj} h_{jj} 
 \Omega\\
& \quad +\frac{2(-1)^{k+1}(n-k-1)(s+2)(s+1)}{(n-1)(n-2)} h \xi_X \Omega,\\
\eqref{X_rho3} \wedge \eqref{D_h_rho3} & = \frac{2(-1)^k(n-k-1)}{(n-1)(n-2)} \sum_{i\neq j, i,j\neq 1} X_{ij}
h_{ij} \Omega.
\end{align*}

Summation of all contributions yields 
\begin{align} 
& X\rho_k \wedge D(h\rho_{n-k}) \nonumber\\
& =(-1)^{k+1}\frac{(s+1)\left[(n-k-2)(n-2)(n+s+1)+2(n-k-1)(n+s)\right]}{(n-1)(n-2)} h \xi_X
\Omega \nonumber \\
& \qquad +\frac{2(-1)^k(n-k-1)}{(n-1)(n-2)} \sum_{i,j=2}^n X_{ij}
h_{ij} \Omega. \label{eq_X_wedge_D}
\end{align}

Using Lemma \ref{lemma_laplacian_hatx}, we write
\begin{align*}
& \sum_{i,j=2}^n X_{ij}h_{ij}\\
& =\sum_{i,j=1}^n X_{ij}h_{ij}-2\sum_{j=1}^n X_{1j}h_{1j}+X_{11}h_{11}\\
& = -\frac12 \Delta_{S^{n-1}} (\hat Xh) +\frac{(n+s)(s+2)}{2} \hat X h-2(s+1) \hat X h+(s+2)(s+1) h
\xi_X\\
& = -\frac12 \Delta_{S^{n-1}} (\hat Xh) +\frac{(n+s)(s+2)-4(s+1)}{2} \hat X h+(s+2)(s+1) h \xi_X.
\end{align*}

Replacing this in  \eqref{eq_X_wedge_D} yields
\begin{align*}
& X\rho_k \wedge D(h\rho_{n-k})\\ 
& =\frac{(-1)^{k+1}(n-k-1)}{(n-1)(n-2)} \Delta_{S^{n-1}}(\hat X h)
\Omega\\
&\quad +\frac{(-1)^{k}(n-k-1)\left[(n+s)(s+2)-4(s+1)\right]}{(n-1)(n-2)}   \hat X h \Omega\\
&\quad +\frac{(-1)^{k+1}(s+1)\left[(n-k-2)(n+s+1)+2(n-k-1)\right]}{n-1}h \xi_X \Omega,
\end{align*}
which completes the proof.
\endproof 

\begin{Proposition}\label{zentraleProp}
Let $f \in \mathcal{H}^n_s$, $h \in \mathcal{H}^n_{s+2}$ and $s \neq 1$. Then 
\begin{displaymath}
 \langle \pd_m(X \mu_{k,f}),\mu_{n-k,h}\rangle=\flag{n}{k} \frac{(-1)^{s+1}(s-1)(n+s+1)(n-k+s)}{\omega_n(n-1)}
\langle
\xi_X { \cdot} f,h\rangle.
\end{displaymath}
\end{Proposition}

\proof
Clearly, for $f\in C^\infty(S^{n-1})$ we have  
\begin{displaymath}
 X(f  \rho_k)=X(f) \rho_k + f X\rho_k.
\end{displaymath}

Since $\hat X f$ is (the restriction of) a homogeneous polynomial of degree $s$, we have $\langle \hat X f,h\rangle=0$.
Using 
\eqref{eq_relation_x_hatx} and \eqref{eq_x_on_product}, we obtain that 
\begin{align*}
 \langle Xf,h\rangle & =s \langle \xi_X { \cdot} f,h\rangle,\\
\langle f,\hat Xh\rangle & = (n+2s+2)\langle \xi_X { \cdot} f,h\rangle.
\end{align*}

Therefore, by \eqref{eq_rho_Dh},  Proposition \ref{prop_action_poincare}, the self-adjointness of the Laplacian and the fact that $f,h$ are harmonic
of degree $s$ resp. $s+2$, we get
\begin{align*}
&\langle \pd_m(X\mu_{k,f}), \mu_{n-k,h}\rangle =  (-1)^{n+k+s} \langle \pd_m(X \mu_{k,f}),\sigma \mu_{n-k,h}\rangle \\
& = (-1)^{k+s} \int X(f\rho_k)
\wedge D(h \rho_{n-k}) \\
& = \flag{n}{k} \frac{(-1)^{k+s}}{\omega_n}\Big[(-1)^{k}  \left\langle
Xf,h-\frac{1}{n-1}\Delta_{S^{n-1}}h\right\rangle + \\
& \qquad   +a \langle f,\Delta_{S^{n-1}} (\hat
Xh) \rangle+b\langle f,\hat Xh\rangle+c\langle f,\xi_X { \cdot} h\rangle\Big]\\
& =  \flag{n}{k} \frac{(-1)^{k+s}}{\omega_n} \Bigg[(-1)^{k+1}\frac{(s+1)(n+s+1)}{n-1} \langle Xf,h\rangle+ \\
& \qquad  +(a
s(s+n-2)+b) \langle
f,\hat Xh\rangle+c\langle \xi_X { \cdot} f,h\rangle\Bigg]\\
& = \flag{n}{k} \frac{(-1)^{s+1}(s-1)(n+s+1)(n-k+s)}{\omega_n(n-1)} \langle \xi_X { \cdot} f,h\rangle,
\end{align*}
which is the asserted relation.
\endproof

\begin{Lemma}\label{nozero}
 For each $s \ge 0$, there are $f \in \mathcal{H}^n_s$, $h \in \mathcal{H}^n_{s+2}$ and $X \in \mathfrak{h}$ such that 
\begin{displaymath}
 \langle \xi_X { \cdot} f,h\rangle \neq 0.
\end{displaymath}
\end{Lemma}

\proof
We argue by contradiction and assume that $\langle \xi_X { \cdot} f,h\rangle = 0$ for all $f \in \mathcal{H}^n_s$, $h \in \mathcal{H}^n_{s+2}$ and $X \in \mathfrak{h}$. 

 Let $f \in \mathcal{H}_s^n$ and $X \in \mathfrak{h}, X \neq 0$. Since $f$ is the restriction of a polynomial of degree $s$ and $\xi_X$ is the restriction of a quadratic polynomial, \eqref{eq_decomposition_sym} implies that $\xi_X \cdot f \in \bigoplus_i \mathcal{H}_{s+2-2i}^n$. 
 
From the assumption, we obtain that in fact  $\xi_X \cdot f \in \bigoplus_i \mathcal{H}_{s-2i}^n$. On the other hand, for any $\tilde f \in \bigoplus_i \mathcal{H}_{s-2-2i}^n$,  we conclude from \eqref{eq_decomposition_sym} that $\xi_X \cdot \tilde f \in \bigoplus_i \mathcal{H}_{s-2i}^n$.

We thus get a  linear map $\bigoplus_i
\mathcal{H}^n_{s-2i} \to \bigoplus_i \mathcal{H}^n_{s-2i}$, $f \mapsto \xi_X { \cdot} f$. Since this map is clearly injective, it
has
to be onto as well.

Since $\xi_X(e_i)=X_{ii}$ and $\tr X=0$, the continuous function $\xi_X$ has a zero $p_0 \in S^{n-1}$.
Choose a non-zero function $g \in \mathcal{H}^n_s$. By using a rotation, we may suppose that $g(p_0) \neq 0$. But then $g$
is not a multiple of $\xi_X$, which is a contradiction. 
\endproof

As explained in Subsection \ref{subsec_outline}, the recursion formula \eqref{eq_quotient_fourier_s} follows  from the lemma and the preceding propositions. It remains to consider the inductive start $s=3$.  

\subsection{Relation to tensor valuations}
\label{subsec_rel_tensor}

We will consider the translation invariant tensor valuations (Minkowski tensors)
\begin{displaymath}
 \Phi_{k,s}(K)= \binom{n-1}{k} \frac{1}{\omega_{n-k+s} s!} \int_{S^{n-1}} y^s \, dS_k(K,y), 
\end{displaymath}
where $0 \leq k \leq n-1$ and $s \geq 0$; moreover $\Phi_{n,0}:=\vol$. We follow the convention from \cite{hug_schneider_schuster_b} and set $\Phi_{k,s}:=0$ whenever
$k \not\in \{0,1,\ldots,n\}$ or $k=n, s \neq 0$. Clearly $\Phi_{k,s}
\in \TVal^{s,\On(n)}$. By \eqref{eq_centroid}, $\Phi_{k,1}=0$ for $k\in\{0,\ldots,n-1\}$, hence we will usually assume that $s \neq 1$.

For $k\in\{0,\ldots,n\}$, $\Phi_{k,0}=\mu_k$ is the $k$-th intrinsic volume. By Weyl's lemma
\cite{weyl_tubes}, for any nonempty $K$ we have
\begin{equation} \label{eq_weyl}
 \Phi_{0,s}(K)=\frac{1}{\omega_{n+s}s!} \int_{S^{n-1}} y^s \, d\vol_{S^{n-1}}(y)=\left\{ \begin{array}{c c}
\frac{1}{\left(\frac{s}{2}\right)!(4\pi)^\frac{s}{2} } Q^\frac{s}{2},& s \text{ even, } \\
0, & s \text{ odd.}\end{array}\right.
\end{equation}

Clearly, each coefficient in the expansion of $\Phi_{k,s}$ as a linear combination of monomials $e_1^{i_1} \cdots e_n^{i_n}$ (with respect to 
some orthonormal basis $e_1,\ldots,e_n$ of $\R^n$) is of the type $\mu_{k,f}$, for some function $f \in \bigoplus_{i \leq s} \mathcal{H}_{ i}^n$. In order to compute the
Alesker-Fourier transform of such a valuation, we first reduce the problem to the case where these functions are
spherical harmonics. 

\begin{Proposition} \label{prop_relation_phi_psi}
\begin{enumerate}
\item For $0 \leq k <n$ and $ s \neq 1$, define 
\begin{displaymath}
\Psi_{k,s}:=\Phi_{k,s}+ \sum_{j=1}^{\lfloor \frac{s}{2}\rfloor}  \frac{(-1)^j \Gamma(\frac{n-k+s}{2})
\Gamma(\frac{n}{2}+s-1-j)}{(4\pi)^j j!\Gamma(\frac{n-k+s}{2}-j)\Gamma(\frac{n}{2}+s-1)}
Q^j \Phi_{k,s-2j}
\end{displaymath}
 and let $\Psi_{n,0}:=\Phi_{n,0}$.  Then $\Psi_{k,s}$ is the component in $\mathcal{H}_s^n$ (equivalently, the trace free part) 
of $\Phi_{k,s}$. In particular, $\Psi_{k,s} \equiv \Phi_{k,s}
\mod Q$.
\item For $0 \leq k <n$ and $ s \neq 1$, $\Phi_{k,s}$ can be written in terms of $\Psi_{k,s'}$ as 
\begin{displaymath}
 \Phi_{k,s}=\Psi_{k,s}+\sum_{j=1}^{\lfloor \frac{s}{2}\rfloor}
\frac{\Gamma\left(\frac{n-k+s}{2}\right)  \Gamma(\frac{n}{2}+s-2j)}{(4\pi)^j j!
\Gamma(\frac{n-k+s}{2}-j)\Gamma(\frac{n}{2}+s-j)} Q^j \Psi_{k,s-2j}.
\end{displaymath}
\end{enumerate}
\end{Proposition}

\proof
The first assertion is an easy computation, which is based on the identity
 \begin{align*}
  &\Delta_{\R^n} (y^{s-2j} \|y\|^{2j})\\
  &\qquad=-(s-2j)(s-2j-1)y^{s-2j-2}\|y\|^{2j} Q+2j(2j-2s-n+2)y^{s-2j}\|y\|^{2j-2}.
\end{align*}

Let us prove the second assertion. First, we remark that except for the case where $j=0,n=2,s=0$ the coefficient in the definition of $\Psi_{k,s}$ is well-defined and equals $1$ for $j=0$. Still this is the reason for starting the summation from $j=1$ and not from $j=0$. However, considering $n$ as a real parameter, we can 
avoid the need to distinguish different cases in the subsequent calculations, since this way the encountered singularity  is a
removable singularity with limit $1$.  

We now compute 
\begin{align*}
\Psi_{k,s} & +\sum_{j=1}^{\lfloor \frac{s}{2}\rfloor}
\frac{\Gamma\left(\frac{n-k+s}{2}\right)  \Gamma(\frac{n}{2}+s-2j)}{(4\pi)^j j!
\Gamma(\frac{n-k+s}{2}-j)\Gamma(\frac{n}{2}+s-j)} Q^j \Psi_{k,s-2j}\\
& = \lim_{\nu \to n} \sum_{j=0}^{\lfloor \frac{s}{2}\rfloor}
\frac{\Gamma\left(\frac{\nu-k+s}{2}\right)  \Gamma(\frac{\nu}{2}+s-2j)}{(4\pi)^j j!
\Gamma(\frac{\nu-k+s}{2}-j)\Gamma(\frac{\nu}{2}+s-j)} Q^j \Psi_{k,s-2j}\allowdisplaybreaks\\
& = \lim_{\nu \to n} \sum_{j=0}^{\lfloor \frac{s}{2}\rfloor}
\frac{\Gamma\left(\frac{\nu-k+s}{2}\right) \Gamma(\frac{\nu}{2}+s-2j)}{(4\pi)^j j!
\Gamma(\frac{\nu-k+s}{2}-j)\Gamma(\frac{\nu}{2}+s-j)} Q^j  \frac{\Gamma(\frac{\nu-k+s-2j}{2})}{\Gamma(\frac{\nu}{2}+s-2j-1)}\cdot\\
& \qquad\qquad \cdot \sum_{l=0}^{\lfloor\frac{s-2j}{2}\rfloor}\frac{(-1)^l
\Gamma(\frac{\nu}{2}+s-2j-l-1)}{(4\pi)^l l!\Gamma(\frac{\nu-k+s-2j}{2}-l)}
Q^l \Phi_{k,s-2j-2l}\allowdisplaybreaks\\
& = \lim_{\nu \to n} \sum_{j=0}^{\lfloor \frac{s}{2}\rfloor}\sum_{l=0}^{\lfloor\frac{s-2j}{2}\rfloor}
(-1)^l\frac{\Gamma\left(\frac{\nu-k+s}{2}\right) (\frac{\nu}{2}+s-2j-1)\Gamma(\frac{\nu}{2}+s-2j-l-1)}{\Gamma(\frac{\nu}{2}+s-j)\Gamma\left(\frac{\nu-k+s}{2}
-j-l\right)}\cdot \\
& \qquad\qquad \cdot \frac{1}{(4\pi)^{j+l}}\frac{1}{j!l!}Q^{j+l}\Phi_{k,s-2(j+l)}\allowdisplaybreaks\\
& = \sum_{m=0}^{\lfloor \frac{s}{2}\rfloor}(-1)^{m}\Bigg\{\lim_{\nu \to n} \frac{\Gamma\left(\frac{\nu-k+s}{2}
\right)}{\Gamma\left(\frac{\nu-k+s}{2}
-m\right)} \sum_{j=0}^{m}
(-1)^{j}\binom{m}{j}\cdot \\
&\qquad\qquad   \cdot \frac{(\frac{\nu}{2}+s-2j-1)\Gamma(\frac{\nu}{2}+s-m-j-1)}{\Gamma(\frac{\nu}{2}+s-j)}
\Bigg\}\frac{1}{(4\pi)^{m}}\frac{1}{m!}Q^{m}\Phi_{k,s-2m}.
\end{align*}

Hence, the assertion follows once the identity 
\begin{displaymath}
 \sum_{j=0}^m (-1)^j\binom{m}{j}\frac{
\left(\frac{\nu}{2}+s-2j-1\right)\Gamma\left(\frac{\nu}{2}+s-m-j-1\right)}{\Gamma\left(\frac{\nu}{2}+s-j\right)}
=\begin{cases} 1, & m =0, \\ 0, & m>0, \end{cases}
\end{displaymath}
is established for $s,m\in\mathbb{N}_0$, $0 \leq m \leq \lfloor\frac{s}{2}\rfloor$ and all real numbers $\nu$ such that the sum is well-defined. 

The identity is clearly true if 
$m=0$. Hence let $m \geq 1$. Let the parameters $s\in\mathbb{N}_0$ and $\nu$ be fixed and define 
\begin{align*}
f(m,j)&:=(-1)^j\binom{m}{j}\frac{ (\frac{\nu}{2}+s-2j-1)\Gamma(\frac{\nu}{2}+s-m-j-1)}{\Gamma(\frac{\nu}{2}+s-j)}
,\\
G(m,j) &:=(-1)^{j+1}\binom{m-1}{j-1} \frac{\Gamma\left(\frac{\nu}{2}+s-m-j\right)}
{\Gamma\left(\frac{\nu}{2}+s-j\right)},
\end{align*}
for $j\in\{0,\ldots,m\}$. Then, in particular, we have $G(m,0)=G(m,m+1)=0$. 
The function $G$  is suggested by Maple's implementation of Zeilberger's algorithm.  By 
elementary 
calculations we confirm that 
\begin{displaymath}
f(m,j)=G(m,j+1)-G(m,j),\qquad j\in\{0,\ldots,m\},
\end{displaymath}
where the cases $j\in \{0,m\}$ are checked separately, and hence
\begin{displaymath}
\sum_{j=0}^mf(m,j)=G(m,m+1)-G(m,0)=0,
\end{displaymath}
which proves the asserted identity.
\endproof

\begin{Corollary} \label{cor_fourier_psi}
For $0<k<n$ and $s \neq 1$, 
 \begin{displaymath}
  \mathbb{F} (\Psi_{k,s})=c_{n,k,s}\frac{\Gamma\left(\frac{k}{2}\right)\Gamma\left(\frac{
s+n-k}{2}\right)}{\Gamma\left(\frac{n-k}{2}\right)\Gamma\left(\frac{s+k}{2}\right)} \Psi_{n-k,s},
 \end{displaymath}
where $c_{n,k,s}$ is the constant from \eqref{eq_schur}.
\end{Corollary}

\proof  It follows from Proposition \ref{prop_relation_phi_psi}  (i) that
\begin{displaymath}
\Psi_{k,s}=\frac{\Gamma\left(\frac{n-k+s}{2}\right)}{\Gamma\left(\frac{n-k}{2}\right)}\frac{1}{\pi^{\frac{s}{2}}s!}\, \mu_{k,f_s}, 
\end{displaymath}
where 
\begin{displaymath}
f_s(y)=1+\sum_{j=1}^{\left\lfloor\frac{s}{2}\right\rfloor}(-1)^j\frac{\Gamma\left(\frac{n}{2}+s-1-j\right)}{4^j j!\Gamma\left(\frac{n}{2}+s-1\right)(s-2j)!}\, Q^jy^{s-2j},\qquad y\in 
S^{n-1}. 
\end{displaymath}
Since the components of $f_s$ are in $\mathcal{H}^n_s$ (see Proposition \ref{prop_relation_phi_psi}  (ii)) and $f_s$ is independent of $k$, we can apply \eqref{eq_schur} to get 
\begin{align*}
\mathbb{F}(\Psi_{k,s})&=\frac{\Gamma\left(\frac{n-k+s}{2}\right)}{\Gamma\left(\frac{n-k}{2}\right)}\frac{1}{\pi^{\frac{s}{2}}s!}\, c_{n,k,s}\mu_{n-k,f_s}\\
&=\frac{\Gamma\left(\frac{n-k+s}{2}\right)}{\Gamma\left(\frac{n-k}{2}\right)}\frac{1}{\pi^{\frac{s}{2}}s!}\, c_{n,k,s}\pi^{\frac{s}{2}}s!
\frac{\Gamma\left(\frac{k}{2}\right)}{\Gamma\left(\frac{k+s}{2}\right)}\Psi_{n-k,s},
\end{align*}
which yields the asserted relation.
\endproof

\begin{Proposition}
Let $\chi^{(1)},\chi^{(2)}$ be the constants from \cite{hug_schneider_schuster_b}. Then, for $1 \leq k \leq n-2$ and $s \neq 1$,   
\begin{align} 
 \frac{c_{n,k,s}}{c_{n,k+1,s}}
 &=
\frac{\Gamma\left(\frac{k+1}{2}\right)\Gamma\left(\frac{n-k}{2}\right)\Gamma\left(\frac{n-k+s-1}{2}\right)\Gamma\left(\frac{n+1}{2}\right)}{\Gamma\left(\frac{n}{2}\right)\Gamma\left(\frac{n-k-1}{2}\right)\Gamma\left(\frac{n-k+s}{2}\right)\Gamma\left(\frac{k+2}{2}\right)}\cdot\nonumber\\
\nonumber\\&\qquad \cdot 
(\chi^{(1)}_{n,k,n-1,s,0}-2\pi s
\chi^{(2)}_{n,k,n-1,s,0}).\label{eq_quotient_fourier_k}
\end{align}
\end{Proposition}

We remark that the constants $\chi^{(1)},\chi^{(2)}$ are rather complicated iterated sums. For this reason, we refrain from writing them down explicitly here.  
\proof
By Theorem 2.3 of \cite{hug_schneider_schuster_b}, we obtain for $1 \leq k \leq n-2$ that
\begin{displaymath} 
 \mu_1 \cdot \Psi_{k,s} \equiv \flag{n}{1} (\chi^{(1)}_{n,k,n-1,s,0}-2\pi s \chi^{(2)}_{n,k,n-1,s,0}) \Psi_{k+1,s} \mod Q.
\end{displaymath}
Both sides of this equation belong to the same isotypical component of $\SO(n)$, which is $\Val_{k+1,s}$. On the other hand, a non-zero multiple of $Q$ of rank $s$ cannot belong to this isotypical component. Therefore
we actually get the equality
\begin{equation} \label{eq_L_Psi}
 \mu_1 \cdot \Psi_{k,s} = \flag{n}{1} (\chi^{(1)}_{n,k,n-1,s,0}-2\pi s \chi^{(2)}_{n,k,n-1,s,0}) \Psi_{k+1,s}.
\end{equation}

Let $\Lambda$ be the derivation operator, which acts on translation invariant valuations $\phi$ by 
\begin{displaymath}
 \Lambda (\phi)(K):=\left.\frac{d}{dt}\right|_{t=0} \phi(K+tB).
\end{displaymath}
Recall from \cite[Corollary 1.7]{bernig_fu06} that $\Lambda (\phi)=2 \mu_{n-1} * \phi$. We
compute  
\begin{align}
 \Lambda(\Phi_{k,s})(K) & =  \binom{n-1}{k} \frac{1}{\omega_{n-k+s} s!} \int y^s
\left.\frac{d}{dt}\right|_{t=0} \,dS_k(K+tB,y) \nonumber\\
& = \binom{n-1}{k} \frac{1}{\omega_{n-k+s} s!} \int y^s k\, 
dS_{k-1}(K,y) \nonumber \\
& = \frac{\omega_{n-k+s+1}(n-k)}{\omega_{n-k+s}} \Phi_{k-1,s}. \label{eq_derivation_minkowski}
\end{align}
By Proposition \ref{prop_relation_phi_psi} (ii), we obtain that
\begin{displaymath}
2\mu_{n-1}*\Psi_{k,s}= \Lambda(\Psi_{k,s})=\frac{\omega_{n-k+s+1}(n-k)}{\omega_{n-k+s}}
\Psi_{k-1,s}.
\end{displaymath}

Applying the Alesker-Fourier transform to \eqref{eq_L_Psi} and using $\mathbb{F}(\mu_1)=\mu_{n-1}$, we arrive at 
\eqref{eq_quotient_fourier_k}.
\endproof

It remains to compute the constants explicitly in the case $s=3$. In this case, the iterated sums in
\cite{hug_schneider_schuster_b} simplify to 
\begin{align*}
 \chi^{(1)}_{n,k,n-1,3,0} & = \frac{(k+1)\Gamma\left(\frac{k+5}{2}\right)\Gamma\left(\frac{n}{2}\right)}{2 \Gamma\left(\frac{n+1}{2}\right)\Gamma\left(\frac{k}{2}+3\right)}, \\
 \chi^{(2)}_{n,k,n-1,3,0} & = \frac{\Gamma\left(\frac{k+3}{2}\right)\Gamma\left(\frac{n}{2}\right)}{8 \pi \Gamma\left(\frac{n+1}{2}\right)\Gamma\left(\frac{k}{2}+3\right)}.
\end{align*}

From \eqref{eq_quotient_fourier_k} it follows that 
\begin{displaymath}
 \frac{c_{n,k,3}}{c_{n,k+1,3}} =\frac{(n-k)(n-k-1)(n-k+1) \Gamma\left(\frac{n-k}{2}\right)^2
\Gamma\left(\frac{k+3}{2}\right)^2}{k(k+2)(k+1)\Gamma\left(\frac{n-k+3}{2}\right)^2\Gamma\left(\frac{k}{2}\right)^2}.
\end{displaymath}

The solution of this equation is 
\begin{displaymath}
 c_{n,k,3}=\epsilon_n \i^3 \frac{\Gamma\left(\frac{3+k}{2}\right)\Gamma\left(\frac{n-k}{2}\right)}{\Gamma\left(\frac{k}{2}\right)\Gamma\left(\frac{
n+3-k}{2}\right)}, \qquad 1\le k\le n-1, 
\end{displaymath}
where $\epsilon_n$ only depends on $n$. 

Together with \eqref{eq_quotient_fourier_s} we find, for odd $s \neq 1$, $n \geq 3$ and $1 \leq k \leq n-1$, that  
\begin{displaymath}
 c_{n,k,s}=\epsilon_n \i^s \frac{\Gamma\left(\frac{n-k}{2}\right)\Gamma\left(\frac{s+k}{2}\right)}{\Gamma\left(\frac{k}{2}\right)\Gamma\left(\frac{
s+n-k}{2}\right)},
\end{displaymath}
hence by Corollary \ref{cor_fourier_psi}
\begin{displaymath}
 \mathbb{F} (\Psi_{k,s})=\epsilon_n \i^s \Psi_{n-k,s}.  
\end{displaymath}

It remains to compute $\epsilon_n$. For this, we use a functorial property of the Alesker-Fourier transform and a version of the
Crofton
formula, where we integrate tensor valuations {\it with respect to  subspaces} $E$ parallel to affine subspaces
$\bar E$. 
In this case the space in which the 
valuation is to be considered is indicated by a superscript. Here we also use the translation invariance of the given
valuations.

\begin{Lemma}\label{Lemma4.13}
If $K\in\mathcal{K}(\R^n)$, then
\begin{displaymath}
\int_{\overline{\Gr}_2(\R^n)} \Phi_{1,3}^E(K \cap \bar E)\, d\bar E = \binom{n}{2}^{-1} \Phi_{n-1,3}(K). 
\end{displaymath} 
\end{Lemma}

\proof
Let $\bar \chi^{(1)}, \bar \chi^{(2)}$ be the constants from \cite{hug_schneider_schuster_b}. From Theorem 2.6 of
\cite{hug_schneider_schuster_b}, it follows that 
\begin{displaymath}
\int_{\overline{\Gr}_2(\R^n)} \Phi_{1,3}^E(K \cap \bar E) \, d\bar E = (\bar \chi^{(1)}_{n,1,2,3,0}-6\pi \bar
\chi^{(2)}_{n,1,2,3,0}) \Phi_{n-1,3}(K). 
\end{displaymath} 
Working out the constant on the right-hand side, we get the result. 
\endproof

\begin{Lemma} \label{lemma_projection_formula}
 Let $\pi_E:\R^n \to E$ be the orthogonal projection to a two-dimensional subspace $E$ of $\R^n$. 
If $K\in\mathcal{K}(\R^n)$, then
\begin{displaymath}
\Phi_{1,3}(K)=\binom{n}{2} \int_{\Gr_2(\R^n)} \Phi_{1,3}^E(\pi_EK) \, dE. 
\end{displaymath}
\end{Lemma}

\proof

A special case of Theorem 4.4.10 in \cite{schneider_book14} is the following. If $f:S^{n-1} \to \R$ is a continuous
function, then 
\begin{displaymath}
 \int_{\Gr_2(\R^n)} \int_{S^{n-1} \cap E} f(y) \, dS_1^E(\pi_EK,y)\, dE=\frac{\omega_2}{\omega_n} \int_{S^{n-1}} f(y)\, 
dS_1(K,y),
\end{displaymath}
where the superscript in $S_1^E$ indicates that the area measure is determined with respect to the subspace $E$. 
Applying this formula to the components of the function $f(y)=y^3$ gives 
\begin{displaymath}
\binom{n}{2} \int_{\Gr_2(\R^n)} \Phi_{1,3}^E(\pi_EK)\, dE = \Phi_{1,3}(K).
\end{displaymath}
\endproof

Let us recall Theorem 6.2.4 from \cite{alesker_fourier}. It was stated in invariant terms, but can be
translated to the following Euclidean version. 
\begin{Proposition} \label{prop_compatibility_pull_backs}
Let $V$ be a Euclidean vector space. If  
$\phi \in \Val_k^{sm,-}(V)$, for $k\in\{1,\ldots,n-1\}$, is given by 
\begin{displaymath}
 \phi(K)=\int_{\Gr_{k+1}(V)} \psi_E(\pi_E K)\, dm(E),
\end{displaymath}
where $E \mapsto \psi_E  \in \Val_k^-(E)$ is a smooth family of odd valuations of degree $k$ and $m$ is a smooth measure, then 
\begin{displaymath}
 \mathbb{F}(\phi)(K)=\int_{\AGr_{k+1}(V)} \mathbb{F} (\psi_E)(K \cap \bar E) \, d\bar
m(\bar E),
\end{displaymath}
where $\bar m$ is the product measure of $m$ and the Lebesgue measure. Note that the Alesker-Fourier transform on the left-hand
side is with respect to $V$, while it is with respect to $E$ on the right-hand side. 
\end{Proposition}

We apply this formula with $k=1$ and $\psi_E=\Phi_{1,3}^E$ to the expression from Lemma
\ref{lemma_projection_formula}. Using the two-dimensional case of Theorem \ref{mainthm_fourier} (which was shown in
Subsection \ref{subsec_twodim}) and Lemma \ref{Lemma4.13}, we obtain  
\begin{align*}
 \mathbb{F} (\Phi_{1,3})& = \binom{n}{2}  \int_{\AGr_2(\R^n)} \mathbb{F}(\Phi_{1,3}^E)(\cdot \cap \bar E)\, d \bar
E\\
& = \binom{n}{2}  \i^3 \int_{\AGr_2(\R^n)} \Phi_{1,3}^E(\cdot \cap \bar E) \, d \bar E\\
& =   \i^3 \Phi_{n-1,3}.
\end{align*}

On the other hand,
$\mathbb{F} (\Phi_{1,3})=\epsilon_n \i^3 \Phi_{n-1,3}$, hence $\epsilon_n=1$ for all $n$. This finishes the
proof of Theorem \ref{mainthm_fourier}.

\section{Algebraic properties of tensor valuations}
\label{sec_alg_prop}

In this section, we give explicit formulas for convolution, Poincar\'e duality, Ales\-ker-Fourier transform and product 
 of tensor valuations.

\subsection{Convolution}

The convolution of translation invariant tensor valuations will be determined explicitly by means of the techniques provided in 
\cite{bernig_fu06}. We will also use some of the calculations involved in Proposition \ref{prop_rumin}, but the Fourier transform is not required for this purpose. 

\begin{Theorem} \label{thm_convolution}
 The convolution product of tensor valuations is given for $k,l \leq n$ with $k+l \geq n$ and $s_1,s_2\neq 1$
by 
\begin{align*}
 \Phi_{k,s_1} * \Phi_{l,s_2} & = \frac{\omega_{s_1+s_2+2n-k-l}}{\omega_{s_1+n-k}\omega_{s_2+n-l}}
\frac{(n-k)(n-l)}{2n-k-l} \cdot \\
& \quad \cdot \binom{2n-k-l}{n-k}  \binom{s_1+s_2}{s_1} \frac{(s_1-1)(s_2-1)}{1-s_1-s_2} 
\Phi_{k+l-n,s_1+s_2}. 
\end{align*}
\end{Theorem}

\proof 
The cases where $k=n$ or $l=n$ can be checked directly. The consistency with the given formula can be seen from
$$
\frac{\omega_{s_1+s_2+2n-k-l}}{\omega_{s_1+n-k}\omega_{s_2+n-l}}
\frac{(n-k)(n-l)}{2n-k-l}=
\frac{\frac{n-k}{2}\Gamma\left(\frac{s_1+n-k}{2}\right)\frac{n-l}{2}\Gamma\left(\frac{s_2+n-l}{2}\right)}
{\frac{2n-k-l}{2}\Gamma\left(\frac{s_1+s_2+2n-k-l}{2}\right)},
$$
since $s_1=0$ if $k=n$ and $s_2=0$ if $l=n$.

Let us recall a formula from \cite{bernig_fu06}. Let $*_V$ be the Hodge star on the space $\Omega^*(V)$ of
differential forms on $V$. Let $\pi_1:SV \to V, \pi_2:SV \to S^{n-1}$ be the projections. 

Let $*_1$ be the
linear operator on $\Omega^*(SV)$ which is uniquely defined by  
\begin{align*}
*_1(\pi_1^* \tau_1 \wedge \pi_2^* \tau_2)&=(-1)^{\binom{n-\deg \tau_1}{2}} \pi_1^* (*_V \tau_1) \wedge
 \pi_2^* \tau_2, \\
 \nonumber  \tau_1 \in \Omega^*(V),&\quad\tau_2 \in \Omega^*(S^{n-1}).
\end{align*} 

Suppose that valuations $\phi_1 \in \Val^{sm}_k,\phi_2 \in \Val_l^{sm}$ with $0 \leq k,l \leq n-1$ and $k+l \geq n$ are
given by differential forms $\gamma_1,\gamma_2$ of bidegree $(k,n-k-1)$ and $(l,n-l-1)$, respectively, on the
sphere
bundle. Then the convolution
product $\phi_1 * \phi_2$ is given by the form 
\begin{equation} \label{eq_bernig_fu}
 \gamma_1\ \hat{*} \ \gamma_2 := *_1^{-1} \left( *_1 \gamma_1 \wedge *_1 D\gamma_2 \right).
\end{equation}

We first study the case $k=l=n-1$. We thus apply the formula to the tensor valuations $\Phi_{n-1,s_1}$ and
$\Phi_{n-1,s_2}$. In this case, 
\begin{displaymath}
 \gamma_1 := \frac{2}{\omega_{s_1+1} s_1!} y^{s_1} \rho_{n-1}, \gamma_2 = \frac{2}{\omega_{s_2+1} s_2!} y^{s_2}
\rho_{n-1},
\end{displaymath}
where $y$ is a vector valued function and the differential forms are tensor valued. 

Let $T:=\sum_{i=1}^n y_i \frac{\partial}{\partial x_i}$ be the Reeb vector field on $S\R^n$. 
It satisfies $i_T \alpha=1$ and $i_T d\alpha=0$. From $*_1 \rho_{n-1} = \frac{(-1)^{n-1}}{2} \alpha$ it follows that 
\begin{equation} \label{eq_convolution_forms}
 \gamma_1\ \hat{*} \ \gamma_2=\frac{2}{\omega_{s_1+1}\omega_{s_2+1} s_1!s_2!}  i_T y^{s_1} D(y^{s_2} \rho_{n-1}). 
\end{equation}

The Rumin differential was computed in Proposition \ref{prop_rumin}, which gives
\begin{equation} \label{eq_rumin_tensorcase}
 D(y^{s_2} \rho_{n-1})=y^{s_2} d\rho_{n-1}+\pi \alpha \wedge d\xi_{n-2,y^{s_2}}.
\end{equation}

We have 
\begin{equation}\label{drhok}
 d\rho_{n-k}=k\frac{\omega_{k+1}}{\omega_k} \alpha \wedge \rho_{n-k-1}. 
\end{equation}

It is easily checked that for $k\in \{0,\ldots,n-2\}$ we have 
\begin{equation} \label{eq_d_xi}
 d\xi_{k,y} \equiv -(n-k-1) y \rho_k \mod \alpha.
\end{equation}

In particular, $d\xi_{n-2,y} \equiv -y \rho_{n-2} \mod \alpha$ and therefore
\begin{align*}
d\xi_{n-2,y^{s_2}} &= d\left(s_2 y^{s_2-1}\xi_{n-2,y}\right)\\
& \equiv s_2(s_2-1)y^{s_2-2} dy \wedge
\xi_{n-2,y}- s_2 y^{s_2} \rho_{n-2} \mod \alpha.
\end{align*}

Similarly, 
\begin{align*}
 d\xi_{n-2,y^{s_1+s_2}} &\equiv (s_1+s_2)(s_1+s_2-1)y^{s_1+s_2-2} dy \wedge
\xi_{n-2,y}\\
&\qquad \qquad -(s_1+s_2) y^{s_1+s_2}
\rho_{n-2} \mod \alpha.
\end{align*}

From these equations it follows that 
\begin{displaymath}
 y^{s_1} d\xi_{n-2,y^{s_2}} \equiv \frac{s_1s_2}{1-s_1-s_2} y^{s_1+s_2} \rho_{n-2},
\end{displaymath}
where we compute modulo the ideal generated by $\alpha, d\alpha$
and exact forms. 

Substituting these formulas into \eqref{eq_convolution_forms} and \eqref{eq_rumin_tensorcase}, 
we obtain that 
\begin{align*}
  \gamma_1\ \hat{*} \ \gamma_2 & \equiv  \frac{2}{\omega_{s_1+1}\omega_{s_2+1} s_1!s_2!} \left(y^{s_1+s_2} i_T
d\rho_{n-1}+\pi y^{s_1}
d\xi_{n-2,y^{s_2}}\right) \\
& \equiv  \frac{2}{\omega_{s_1+1}\omega_{s_2+1} s_1!s_2!} \left(\pi y^{s_1+s_2} \rho_{n-2}+\pi
\frac{s_1s_2}{1-s_1-s_2}
y^{s_1+s_2} \rho_{n-2}\right)\\
& \equiv  \frac{2 \pi}{\omega_{s_1+1}\omega_{s_2+1} s_1!s_2!}  \frac{(s_1-1)(s_2-1)}{1-s_1-s_2} y^{s_1+s_2}\rho_{n-2},
\end{align*}
which translates to 
\begin{align*}
 &\Phi_{n-1,s_1} * \Phi_{n-1,s_2}\\ 
 & = \frac{2\pi}{\omega_{s_1+1}\omega_{s_2+1} s_1!s_2!} 
\frac{(s_1-1)(s_2-1)}{1-s_1-s_2}
\frac{\omega_{s_1+s_2+2}(s_1+s_2)!}{2 \pi} \Phi_{n-2,s_1+s_2}\\
&  = \binom{s_1+s_2}{s_1} \frac{\omega_{s_1+s_2+2}}{\omega_{s_1+1}\omega_{s_2+1}} 
\frac{(s_1-1)(s_2-1)}{1-s_1-s_2}
 \Phi_{n-2,s_1+s_2}.
\end{align*}

The general case can be derived from this as follows. Repeated application of \eqref{eq_derivation_minkowski}
shows that, 
for $0 \leq a \leq n-2$, we have
\begin{displaymath}
 \Lambda^a (\Phi_{n-1,s})=\frac{\omega_{s+a+1}}{\omega_{s+1}} a! \Phi_{n-a-1,s}.
\end{displaymath}

Since $\Lambda$ acts as convolution by $2\mu_{n-1}$, compare \cite{bernig_fu06}, we obtain 
\begin{align*}
& \Phi_{k,s_1} * \Phi_{l,s_2} \\
& =  \frac{\omega_{s_1+1}}{\omega_{s_1+n-k}(n-k-1)!}  \Lambda^{n-k-1}\Phi_{n-1,s_1} * 
\frac{\omega_{s_2+1}}{\omega_{s_2+n-l}(n-l-1)!}  \Lambda^{n-l-1}\Phi_{n-1,s_2}\\
& = \frac{\omega_{s_1+1}}{\omega_{s_1+n-k}(n-k-1)!} \frac{\omega_{s_2+1}}{\omega_{s_2+n-l}(n-l-1)!}
\binom{s_1+s_2}{s_1} \cdot \\
& \quad \cdot \frac{\omega_{s_1+s_2+2}}{\omega_{s_1+1} \omega_{s_2+1}} 
\frac{(s_1-1)(s_2-1)}{1-s_1-s_2} \Lambda^{2n-k-l-2} \Phi_{n-2,s_1+s_2}\\
& = \frac{\omega_{s_1+1}}{\omega_{s_1+n-k}(n-k-1)!} \frac{\omega_{s_2+1}}{\omega_{s_2+n-l}(n-l-1)!} \binom{s_1+s_2}{s_1}\frac{\omega_{s_1+s_2+2}}{\omega_{s_1+1} \omega_{s_2+1}}
\cdot\\
& \quad \cdot  
\frac{(s_1-1)(s_2-1)}{1-s_1-s_2} \frac{\omega_{s_1+s_2+2n-k-l}(2n-k-l-1)!}{\omega_{s_1+s_2+2}} \Phi_{k+l-n,s_1+s_2}\\
& = \frac{\omega_{s_1+s_2+2n-k-l}}{\omega_{s_1+n-k}\omega_{s_2+n-l}} \frac{(n-k)(n-l)}{2n-k-l} \binom{2n-k-l}{n-k} 
\binom{s_1+s_2}{s_1} \cdot\\
& \quad \cdot 
\frac{(s_1-1)(s_2-1)}{1-s_1-s_2}  \Phi_{k+l-n,s_1+s_2},
\end{align*}
which yields the assertion of the theorem.
\endproof

\subsection{Poincar\'e duality}

\begin{Proposition} \label{prop_poincare_tensorcase}
For $0 \leq k \leq n$ and $ s \neq 1$, the Poincar\'e pairing with respect to  the convolution is given by
\begin{displaymath}
 \langle \pd_c^s(\Phi_{k,s}),\Phi_{n-k,s}\rangle= 
\frac{1-s}{\pi^s s!^2}\binom{n}{k}  \frac{k(n-k)}{4}
\frac{\Gamma\left(\frac{k+s}{2}\right)\Gamma\left(\frac{n-k+s}{2}\right)}{\Gamma\left(\frac{n}{2}+1\right)}.
\end{displaymath}
\end{Proposition}

\proof
In the cases $k=0,n$ (with $s=0$) the formula is easily confirmed. The right-hand  side is equal to one, if properly
interpreted.

Let us assume $1 \leq k \leq n-1$. Recall that $\Phi_{k,s}$ is represented 
by $\gamma_1:=\frac{\omega_{n-k}}{\omega_{n-k+s}s!} y^s \rho_k$ and $\Phi_{n-k,s}$ is 
represented by $\gamma_2:=\frac{\omega_k}{\omega_{k+s}s!} y^s \rho_{n-k}$.  

By Proposition \ref{prop_rumin}, \eqref{eq_d_xi} and the subsequent considerations, and by \eqref{drhok}, we get

\begin{align*}
D(y^s \rho_{n-k}) & =y^s d\rho_{n-k}+\frac{\omega_{k+1}}{\omega_k} \alpha \wedge d\xi_{n-k-1,y^s} \\
& = y^s d\rho_{n-k}+\frac{\omega_{k+1}}{\omega_k} \alpha \wedge \left(s(s-1)y^{s-2} dy \wedge
\xi_{n-k-1,y}- s k y^s \rho_{n-k-1}\right)\\
& = (1-s)y^s D\rho_{n-k}+\frac{\omega_{k+1}}{\omega_k} \alpha \wedge s(s-1)y^{s-2} dy \wedge
\xi_{n-k-1,y}.
\end{align*}

Since $\sum y_i^2=1$, we obtain that $\contr(y^s,y^s)=1$ and $\contr(y^s,dy \wedge \tau)=0$ for all $\mathrm{Sym}^{s-1}V$-valued forms $\tau$. Therefore the contraction of $y^s$ with the second term vanishes and we obtain 
\begin{align*}
 \contr(\gamma_1,D\gamma_2) & =\frac{\omega_{n-k}}{\omega_{n-k+s}s!} \frac{ \omega_k}{\omega_{k+s}s!} \contr(y^s \rho_k,D(y^s \rho_{n-k}))\\
& =  \frac{\omega_{n-k}}{\omega_{n-k+s}s!} \frac{ \omega_k}{\omega_{k+s}s!}(1-s)  \rho_k \wedge D\rho_{n-k}\\
& =  (-1)^k \flag{n}{k} \frac{1}{\omega_n} \frac{\omega_{n-k}}{\omega_{n-k+s}s!} \frac{ \omega_k}{\omega_{k+s}s!}(1-s)
\vol_{S\R^n},
\end{align*}
where we used \eqref{eq_rho_Dh} in the last line. The statement now follows from Equation (43) in \cite{wannerer_area_measures} (note the different meaning of the constant $\omega_n$ in \cite{wannerer_area_measures}). 
\endproof

\begin{Corollary} \label{cor_poincare_prod}
For $0 \leq k \leq n$ and $s\neq 1$, the Poincar\'e pairing with respect to  the product is  given by
\begin{displaymath}
 \langle \pd_m^s(\Phi_{k,s}),\Phi_{n-k,s}\rangle=  (-1)^s \frac{1-s}{\pi^s s!^2}\binom{n}{k} 
\frac{k(n-k)}{4}
\frac{\Gamma\left(\frac{k+s}{2}\right)\Gamma\left(\frac{n-k+s}{2}\right)}{\Gamma\left(\frac{n}{2}+1\right)}.
\end{displaymath}
\end{Corollary}

\proof
Immediate from Proposition \ref{prop_poincare_tensorcase} and Lemma \ref{lemma_relation_pds}.  
\endproof

We remark that the Poincar\'e duality of multiples of the basic invariants $\Phi_{k,s}$ with powers of $Q$ can be easily written down
by using \eqref{eq_adjoint_Q} and  
\begin{equation} \label{eq_traces}
 \tr \Phi_{k,s}=\frac{n-k+s-2}{2\pi s(s-1)}\Phi_{k,s-2}, \quad s \geq 2. 
\end{equation}

\subsection{Alesker-Fourier transform}
\label{subsec_fourier}

As a straightforward consequence of the results from Section \ref{sec_fourier}, in particular of Corollary \ref{cor_fourier_psi} and the fact 
that $\epsilon_n=1$ by the arguments at the end of Section \ref{sec_fourier}, we obtain the following theorem.

\begin{Theorem} \label{thm_fourier_tensorcase}
For $0 \leq k \leq n$ and $s \neq 1$, the Alesker-Fourier transform of tensor valuations is given by
\begin{align*}
\mathbb{F} (\Psi_{k,s}) & =\i^s \Psi_{n-k,s},\\ 
\mathbb{F} (\Phi_{k,s}) & =\i^s \sum_{j=0}^{\lfloor \frac{s}{2}\rfloor} \frac{(-1)^j}{(4\pi)^j j!}  Q^j
\Phi_{n-k,s-2j}.
\end{align*}
\end{Theorem}

\proof
Note first that $\Psi_{0,s}=0$ and $\Psi_{n,s}=0$ if $s>0$, while $\Psi_{0,0}=\chi$ and $\Psi_{n,0}=\vol_n$. The cases $k=0$ and $k=n$ are thus directly confirmed, so let us assume that $1 \leq k \leq n-1$.  

The first equation follows from Corollary \ref{cor_fourier_psi} and \eqref{eq_fourier_schur}, since we have seen
in the previous section that $\epsilon_n=1$. For the second equation,
we use Proposition \ref{prop_relation_phi_psi} and compute
\begin{align*}
 \mathbb{F}(\Phi_{k,s}) & = \Gamma\left(\frac{n-k+s}{2}\right) \sum_{j=0}^{\lfloor \frac{s}{2}\rfloor}
\frac{\Gamma(\frac{n}{2}+s-2j)}{(4\pi)^j j!
\Gamma(\frac{n-k+s}{2}-j)\Gamma(\frac{n}{2}+s-j)} Q^j \i^{s-2j}\Psi_{n-k,s-2j}\\
& = \Gamma\left(\frac{n-k+s}{2}\right) \lim_{\nu \to n} \sum_{j=0}^{\lfloor \frac{s}{2}\rfloor}
\frac{\Gamma(\frac{\nu}{2}+s-2j)}{(4\pi)^j j!
\Gamma(\frac{\nu-k+s}{2}-j)\Gamma(\frac{\nu}{2}+s-j)} Q^j \i^{s-2j}\cdot \\ 
& \quad \cdot  \frac{\Gamma(\frac{k+s-2j}{2})}{\Gamma(\frac{\nu}{2}+s-2j-1)} \sum_{i=0}^{\lfloor \frac{s}{2}\rfloor-j}
\frac{(-1)^i
\Gamma(\frac{\nu}{2}+s-2j-i-1)}{(4\pi)^i i!\Gamma(\frac{k+s}{2}-i-j)}
Q^i \Phi_{k,s-2j-2i}\\
& = \i^s \Gamma\left(\frac{n-k+s}{2}\right) \lim_{\nu \to n} \sum_{m=0}^{\lfloor \frac{s}{2}\rfloor} 
\frac{(-1)^m}{(4\pi)^m m! \Gamma(\frac{k+s}{2}-m)} \\
& \quad \sum_{j=0}^m  \binom{m}{j}
\frac{(\frac{\nu}{2}+s-2j-1)\Gamma(\frac{k+s}{2}-j)\Gamma(\frac{\nu}{2}+s-j-m-1)}{
\Gamma(\frac{\nu-k+s}{2}-j)\Gamma(\frac{\nu}{2}+s-j)}  
Q^m \Phi_{k,s-2m}\\
& = \i^s  \sum_{m=0}^{\lfloor \frac{s}{2}\rfloor} 
\frac{(-1)^m}{(4\pi)^m m!} 
Q^m \Phi_{k,s-2m},
\end{align*}
where we used that
\begin{displaymath}
\sum_{j=0}^m  \binom{m}{j}
\frac{(\frac{\nu}{2}+s-2j-1)\Gamma(\frac{k+s}{2}-j)\Gamma(\frac{\nu}{2}+s-j-m-1)}{
\Gamma(\frac{\nu-k+s}{2}-j)\Gamma(\frac{\nu}{2}+s-j)}
=\frac{\Gamma(\frac{k+s}{2}-m)}{\Gamma(\frac{\nu-k+s}{2})}
\end{displaymath}
for all real $\nu$ such that the sum is well-defined. 

To verify this, we define for $\nu$ and $s,m,j,k$ as above
\begin{displaymath}
f(m,j):=\binom{m}{j}
\frac{(\frac{\nu}{2}+s-2j-1)\Gamma\left(\frac{k+s}{2}-j\right) \Gamma\left(\frac{\nu}{2}+s-j-m-1\right)}{
\Gamma\left(\frac{\nu-k+s}{2}-j\right)\Gamma\left(\frac{\nu}{2}+s-j\right)}.
\end{displaymath}

Then Maple's implementation of Zeilberger's algorithm provides the auxiliary function 
\begin{displaymath}
\displaystyle G(m,j) := -2\binom{m}{j-1}\frac {   \Gamma\left(
\frac{\nu}{2}+s-m-j-1\right) 
 \Gamma\left(\frac{k+s}{2}-j+1\right) }{ \Gamma\left(\frac{\nu-k+s}{2}-j\right)  
 \Gamma\left(\frac{\nu}{2}+s-j\right)} 
\end{displaymath}
which satisfies $G(m,0)=G(m,m+2)=0$ and
\begin{displaymath}
-2f(m,j)+(k+s-2(m+1))f(m+1,j)=G(m,j+1)-G(m,j) 
\end{displaymath}
for $j\in\{0,\ldots,m+1\}$ and $m+1\le\lfloor \frac{s}{2}\rfloor$. Here the case $j=m+1$ has to be treated separately. 
Summation over $j$ thus yields for $f(m):=\sum_{j=0}^mf(m,j)$ the recursion
\begin{displaymath}
f(m+1)=\frac{1}{\frac{k+s}{2}-(m+1)}f(m)
\end{displaymath}
for $m\in\mathbb{N}_0$. Moreover, we have 
$f(0)=\frac{\Gamma\left(\frac{k+s}{2}\right)}{\Gamma\left(\frac{\nu-k+s}{2}\right)}$. Thus we get 
\begin{displaymath}
f(m)=\frac{\Gamma(\frac{k+s}{2}-m)}{\Gamma(\frac{\nu-k+s}{2})},
\end{displaymath}
which is the asserted identity.
\endproof

\subsection{Product}

Using the explicit formulas for  the convolution and  the Alesker-Fourier transform from Theorems \ref{thm_convolution} and
\ref{thm_fourier_tensorcase}, we can finally compute the product of tensor valuations. 

\begin{Theorem} \label{thm_product_tensorcase}
For $0 \leq k,l$ with $k+l \leq  n$ and $s_1,s_2\neq 1$, the product of tensor valuations is given  by 
\begin{align*}
 \Phi_{k,s_1} \cdot \Phi_{l,s_2} & = \frac{kl}{k+l} \binom{k+l}{k}  \sum_{\substack{a=0\\2a \neq s_1+s_2-1}}^{\lfloor
\frac{s_1+s_2}{2}\rfloor}
\frac{1}{(4\pi)^a a!} \Bigg( \sum_{m=0}^a \sum_{i=\max \left\{0, m-\left\lfloor\frac{s_2}{2}\right\rfloor\right\}}^{\min
\left\{m,\left\lfloor \frac{s_1}{2}\right\rfloor\right\}} \\
& \quad  (-1)^{a-m}
 \binom{a}{m} \binom{m}{i} 
\frac{\omega_{s_1+s_2-2m+k+l}}{\omega_{s_1-2i+k}\omega_{s_2-2m+2i+l}} \binom{s_1+s_2-2m}{s_1-2i} 
 \cdot \\
& \quad \cdot \frac{(s_1-2i-1)(s_2-2m+2i-1)}{1-s_1-s_2+2m} \Bigg) Q^a
\Phi_{k+l,s_1+s_2-2a}.
\end{align*}
\end{Theorem}

\proof
The cases where $k=0$ or $l=0$ can be checked directly. Hence we restrict ourselves to $k,l\ge 1$, in the following. 
For $i=1,2$, set $\delta_i=0$ if $s_i$ is even and $\delta_i=\frac32$ if $s_i$ is odd. By Theorem \ref{thm_fourier_tensorcase}
\begin{align*}
 \Phi_{k,s_1} \cdot \Phi_{l,s_2} & = \mathbb{F}^{-1} \circ \mathbb{F}( \Phi_{k,s_1} \cdot \Phi_{l,s_2} ) \\
& = \mathbb{F}^{-1} \left(\mathbb{F}(\Phi_{k,s_1}) * \mathbb{F} (\Phi_{l,s_2})\right)\\
& = \i^{s_1+s_2} \mathbb{F}^{-1} \left(\sum_{i=0}^{\frac{s_1}{2}-\delta_1}
\sum_{j=0}^{
\frac{s_2}{2}-\delta_2} \frac{(-1)^{i+j}}{(4\pi)^{i+j} i! j!} Q^i \Phi_{n-k,s_1-2i} * Q^j \Phi_{n-l,s_2-2j}\right).
\end{align*}

Then Theorem \ref{thm_convolution} and  the change of variables $m:=i+j$ yield
\begin{align*}
& \Phi_{k,s_1} \cdot \Phi_{l,s_2}\\ 
& =\i^{s_1+s_2} \frac{kl}{k+l} \binom{k+l}{k}  \mathbb{F}^{-1}
\Bigg(\sum_{i,j} Q^{i+j}
\frac{(-1)^{i+j}}{(4\pi)^{i+j} i! j!} 
\frac{\omega_{s_1-2i+s_2-2j+k+l}}{\omega_{s_1-2i+k}\omega_{s_2-2j+l}}
 \cdot \\
& \quad \cdot  \binom{s_1+s_2-2i-2j}{s_1-2i} \frac{(s_1-2i-1)(s_2-2j-1)}{1-s_1-s_2+2i+2j}
\Phi_{n-k-l,s_1+s_2-2i-2j}\Bigg) \\
& = \i^{s_1+s_2} \frac{kl}{k+l} \binom{k+l}{k}  \sum_{m=0}^{\frac{s_1+s_2}{2}-\delta_1-\delta_2} Q^m
\frac{(-1)^m}{(4\pi)^m m!}  \cdot \\ &\quad\cdot \Bigg(\sum_{i=\max \left\{0,m-\frac{s_2}{2}+\delta_2\right\}}^{\min
\left\{m,\frac{s_1}{2}-\delta_1\right\}} \binom{m}{i} 
\frac{\omega_{s_1+s_2-2m+k+l}}{\omega_{s_1-2i+k}\omega_{s_2-2m+2i+l}}
 \cdot \\
& \quad \cdot  \binom{s_1+s_2-2m}{s_1-2i} \frac{(s_1-2i-1)(s_2-2m+2i-1)}{1-s_1-s_2+2m} \Bigg) 
\mathbb{F}^{-1}\left( \Phi_{n-k-l,s_1+s_2-2m}\right),
\end{align*}
where $s_1+s_2-2m\neq 1$. 

Now we apply once again Theorem \ref{thm_fourier_tensorcase} and obtain 
\begin{align*}
\Phi_{k,s_1} \cdot \Phi_{l,s_2} & = \frac{kl}{k+l} \binom{k+l}{k}  \sum_{m=0}^{\frac{s_1+s_2}{2}-\delta_1-\delta_2} Q^m
\frac{1}{(4\pi)^m m!}  \cdot\\&\quad\cdot \Bigg(\sum_{i=\max \left\{0,m-\frac{s_2}{2}+\delta_2\right\}}^{\min
\left\{m,\frac{s_1}{2}-\delta_1\right\}}  \binom{m}{i} 
\frac{\omega_{s_1+s_2-2m+k+l}}{\omega_{s_1-2i+k}\omega_{s_2-2m+2i+l}}
 \cdot \\
& \quad \cdot  \binom{s_1+s_2-2m}{s_1-2i} \frac{(s_1-2i-1)(s_2-2m+2i-1)}{1-s_1-s_2+2m} \Bigg)\cdot \\
& \quad \cdot \sum_{\substack{j=0\\2j \neq s_1+s_2-2m-1}}^{\left\lfloor\frac{s_1+s_2-2m}{2}\right\rfloor} \frac{(-1)^j}{(4\pi)^j j!} Q^j \Phi_{k+l,s_1+s_2-2m-2j}.
\end{align*}
If $m>\frac{s_1+s_2}{2}-\delta_1-\delta_2$, then the inner sum is empty. Hence we can restrict  $m$ to the range from $0$ to $\left\lfloor\frac{s_1+s_2}{2}\right\rfloor$ without changing the sum. Similarly, if $2i=s_1-1$ or $2i=2m-s_2+1$, then the corresponding term in the inner sum vanishes. Hence we can sum over all $\max\left\{0,m-\left\lfloor \frac{s_2}{2}\right\rfloor\right\} \leq i \leq \min\left\{m,\left\lfloor \frac{s_1}{2}\right\rfloor\right\}$. Using the change of variables $a:=m+j$, we get  
\begin{align*}
& \Phi_{k,s_1} \cdot \Phi_{l,s_2}\\
& = \frac{kl}{k+l} \binom{k+l}{k}  \sum_{\substack{a=0\\2a \neq s_1+s_2-1}}^{\left\lfloor \frac{s_1+s_2}{2}\right\rfloor}
\frac{1}{(4\pi)^a a!} \Bigg( \sum_{m=0}^a \sum_{i=\max\left\{0,m-\left\lfloor \frac{s_2}{2}\right\rfloor\right\}}^{\min\left\{m,\left\lfloor \frac{s_1}{2}\right\rfloor \right\}} \\
& \quad  (-1)^{a-m}
 \binom{a}{m} \binom{m}{i} 
\frac{\omega_{s_1+s_2-2m+k+l}}{\omega_{s_1-2i+k}\omega_{s_2-2m+2i+l}}
 \cdot \\
& \quad \cdot  \binom{s_1+s_2-2m}{s_1-2i} \frac{(s_1-2i-1)(s_2-2m+2i-1)}{1-s_1-s_2+2m} \Bigg) Q^a
\Phi_{k+l,s_1+s_2-2a},
\end{align*}
which proves the result.
\endproof

There seems to be no closed formula for the inner double sum, except for small values of $s_1,s_2$. For
further use, we give  some explicit formulas in the case of tensors of small ranks. 

\begin{Corollary} \label{cor_product_tensor_small_rank}
\begin{enumerate}
 \item 
If $0 \leq k,l$ with $k+l \leq n$, then 
\begin{align*}
 \Phi_{k,2} \cdot \Phi_{l,2} & = -\frac{\Gamma\left(\frac{k+l+1}{2}\right)}{\pi^{3/2} (k+l+2)(k+l)
\Gamma\left(\frac{k+1}{2}\right)
\Gamma\left(\frac{l+1}{2}\right)}\cdot\\
& \quad \cdot\left(2\pi^2 kl \Phi_{k+l,4}-\frac{k^2+l^2 +4k l+2k+2 l}{4} \pi Q
\Phi_{k+l,2}+\frac{kl}{16} Q^2 \Phi_{k+l,0} \right).
\end{align*}
\item If $1 \leq k,l$ with $k+l \leq n$, then
\begin{align*}
\Phi_{k,3} \cdot \Phi_{l,3} & =
\frac{(k+1)(l+1)\Gamma\left(\frac{k+l+1}{2}\right)}{\pi^\frac52
(k+l+4)(k+l+2)(k+l)\Gamma\left(\frac{k}{2}\right)\Gamma\left(\frac{l}{2}\right)}\cdot\\
& \quad \cdot\left(-32 \Phi_{k+l,6}\pi^3+8Q\Phi_{k+l,4}\pi^2-Q^2 \Phi_{k+l,2}\pi+\frac{1}{12}Q^3\Phi_{k+l,0}\right).
\end{align*}
\item If  $1 \leq k, 0 \leq l$, then
\begin{align*}
 \Phi_{k,3} \cdot \Phi_{l,2} & = \frac{(k+1)\Gamma\left(\frac{k+l}{2}\right)}{4\sqrt{\pi} 
 (k+l+3)(k+l+1)  \Gamma\left(\frac{k}{2}\right)\cdot \Gamma\left(\frac{l+1}{2} \right)}\cdot\\
& \quad \cdot(-20 \pi l \Phi_{k+l,5}+(k+6l+3)Q \Phi_{k+l,3}).
\end{align*}
\end{enumerate}
\end{Corollary}

\proof
The case $1 \leq k,l$ of all three equations follows by simplifying the formula from Theorem \ref{thm_product_tensorcase}. The remaining cases follow by a direct computation based on \eqref{eq_weyl}. 
\endproof

\section{Explicit kinematic formulas}
\label{sec_explicit_kinforms}

In this section, we translate the algebraic results from the last two sections into geometric statements about
kinematic formulas. 
 
\subsection{Crofton formulas}
\label{subsec_croftonformulas}

The special case of Theorem \ref{thm_product_tensorcase} where one of the two tensor valuations is an intrinsic volume
can be translated to Crofton-type formulas. 

\proof[Proof of Theorem \ref{mainthm_crofton}]
By Theorem \ref{thm_product_tensorcase} we get
\begin{align*}
 \Phi_{k,s} \cdot \mu_l & = \Phi_{k,s} \cdot \Phi_{l,0} \\
& = \binom{k+l}{k} \frac{kl}{k+l} \sum_{a=0,2a \neq s-1}^{\left\lfloor \frac{s}{2}\right\rfloor} 
 \frac{1}{(4\pi)^a a!}\cdot \\
 &\qquad \cdot \sum_{m=0}^a (-1)^{a-m} \binom{a}{m} \frac{\omega_{s-2m+k+l}}{\omega_{s-2m+k}\omega_l} Q^a \Phi_{k+l,s-2a}\\
& = \binom{k+l}{k} \frac{kl}{k+l} \sum_{a=0,2a \neq s-1}^{\left\lfloor \frac{s}{2}\right\rfloor} 
 \frac{1}{(4\pi)^a a!} \frac{\Gamma\left(\frac{l+2a}{2}\right)\Gamma\left(\frac{s+k-2a}{2}\right)}{2
\Gamma\left(\frac{s+k+l}{2}\right)} Q^a \Phi_{k+l,s-2a}.
\end{align*}

Here we used the equation 
\begin{displaymath}
 \sum_{m=0}^a (-1)^{a-m} \binom{a}{m}
\frac{\omega_{s-2m+k+l}}{\omega_{s-2m+k}\omega_l}=\frac{\Gamma\left(\frac{l+2a}{2}\right)\Gamma\left(\frac{s+k-2a}{2}
\right)}{2 \Gamma\left(\frac{s+k+l}{2}\right)}.
\end{displaymath}

To verify this, we show that 
\begin{displaymath}
f(a):=\sum_{m=0}^a(-1)^m\binom{a}{m}\frac{\Gamma\left(\frac{s+k}{2}-m\right)\Gamma\left(\frac{l}{2}\right)}{\Gamma\left(\frac{s+k+l}{2}-m\right)}=(-1)^a\frac{\Gamma\left(\frac{l}{2}+a\right)\Gamma\left(\frac{s+k}{2}-a\right)}{\Gamma\left(\frac{s+k+l}{2}\right)},
\end{displaymath}
for $a\in\{0,\ldots,\lfloor\frac{s}{2}\rfloor\}$ with $2a\neq s-1$ and $k,l\ge 1$.  

A continuity argument shows that it
is sufficient to prove the assertion for real $s>2a$ such that $2a\neq s-1$. To verify the equation for these cases, we
define
\begin{align*}
f(a,m)&:=(-1)^m\binom{a}{m}\frac{\Gamma\left(\frac{s+k}{2}-m\right)\Gamma\left(\frac{l}{2}\right)}{\Gamma\left(\frac{s+k+l}{2}-m\right)},\\
G(a,m)&:=(-1)^m 2 \binom{a}{m-1}\frac{\Gamma\left(\frac{s+k}{2}-m+1\right)\Gamma\left(\frac{l}{2}\right)}{\Gamma\left(\frac{s+k+l}{2}-m\right)},
\end{align*}
for $m=0,\ldots,a+1$  with $a+1\le \lfloor\frac{s}{2}\rfloor$. In addition, we put $ G(a,a+2):=0$. Note that $\frac{s+k+l}{2}>a+1$, since
$s>2a$. Therefore, we have $f(a,a+1)=0$. 
As before we use that $\binom{a}{r}=0$ if $r<0$ or $r>a$. One confirms by a direct calculation that 
\begin{displaymath}
-(2a+l)f(a,m)+(-s-k+2+2a)f(a+1,m)=G(a,m+1)-G(a,m),
\end{displaymath}
for $m=0,\ldots,a+1$, where the  case $m=a+1$ has to be checked separately. Summation over
$m\in\{0,\ldots,a+1\}$ yields that
\begin{displaymath}
-\left(\frac{l}{2}+a\right)f(a)+\left(1+a-\frac{s+k}{2}\right)f(a+1)=0,
\end{displaymath}
that is, 
\begin{displaymath}
f(a+1)=(-1)\frac{\frac{l}{2}+a}{\frac{s+k}{2}-a-1}f(a). 
\end{displaymath}
Note that $\frac{s+k}{2} > a+1$. Using this relation recursively and 
\begin{displaymath}
f(0)=\frac{\Gamma\left(\frac{s+k}{2}\right)\Gamma\left(\frac{l}{2}\right)}{\Gamma\left(\frac{s+k+l}{2}\right)},
\end{displaymath}
the assertion follows.

Now the statement of the theorem is a direct consequence of this formula and the definition of the product of valuations (see \cite[Section 3.5]{bernig_aig10}). 
\endproof

{ We remark that another approach (which indeed leads to the same formula) is to use the recent
results from \cite{goodey_weil14} on $k$-th section bodies, compare also \cite{schuster_wannerer14}.}  

Looking at the trace free-part of the formula in Theorem \ref{mainthm_crofton} (or by a direct argument similar to the one 
for Theorem \ref{mainthm_crofton}), we obtain a Crofton-type formula for the tensor valuation $\Psi_{k,s}$.

\begin{Corollary}[Crofton formula in the $\Psi$-basis] \label{cor_crofton_psi} 
If  $k,l \geq 0$ and $ k+l \leq n$, then
\begin{displaymath}
 \int_{\AGr_{n-l}(\R^n)} \Psi_{k,s}(K \cap \bar E) \, d\bar E=\frac{\omega_{s+k+l}}{\omega_{s+k}\omega_l}
\binom{k+l}{k} \frac{kl}{k+l}  \flag{n}{l}^{-1} \Psi_{k+l,s}(K).
\end{displaymath}
\end{Corollary}

\subsection{Additive formulas}

\proof[Proof of Theorem \ref{thm_additive_kf}]
Let $\Area$ denote the space of smooth area measures \cite{wannerer_area_measures, wannerer_unitary_module}. Roughly
speaking, an area measure is a valuation with values in the space of signed measures on the unit sphere. A smooth area
measure is one that can be represented by integration over the normal cycle of some differential form. 

If $G$ is a subgroup of $\On(n)$ acting transitively on the unit sphere, then $\Area^G$ is
finite-dimensional. Wannerer \cite{wannerer_area_measures} has shown the existence of a kinematic formula 
for area measures of the following type. Let
$\Phi \in \Area^G$, and let $\Psi_1,\ldots,\Psi_m$ be a basis of $\Area^G$. Then there exist constants
$c_{kl}^\Phi$,
depending only on $\Phi$, such that for compact convex bodies $K,L$ and Borel subsets $U_1,U_2 \subset S^{n-1}$ we have
\begin{displaymath}
 \int_G \Phi(K + gL, U_1 \cap g U_2)\, dg=\sum_{k,l=1}^m c_{kl}^\Phi \Psi_k(K,U_1) \Psi_l(L,U_2).
\end{displaymath}
The corresponding map $A:\Area^G \to
\Area^G \otimes \Area^G$ is called {\it additive kinematic  operator for $G$-equivariant area measures}. 

The map
\begin{displaymath}
 M^s:\Area^G \to \TVal^{s,G},\quad  M^s(\Phi)(K):=\int_{S^{n-1}} y^s \,d\Phi(K,y),
\end{displaymath}
is called  the {\it $s$-th moment map}.

The diagram 
\begin{displaymath}
 \xymatrix{\Area^G \ar[r]^<<<<<<{A} \ar[d]_{M^{s_1+s_2}} & \Area^G \otimes
\Area^G
\ar[d]_{M^{s_1} \otimes M^{s_2}} \\
\TVal^{s_1+s_2,G} \ar[r]^<<<<<{a^G_{s_1,s_2}} & \TVal^{s_1,G} \otimes
\TVal^{s_2,G}}
\end{displaymath}
commutes  \cite[Prop. 6]{wannerer_area_measures}.

In the case $G=\On(n)$, $\Area^{\On(n)}$ is spanned by the area measures $S_i$, where $ i=0,\ldots,n-1$.
By
definition of the moment map, 
\begin{displaymath}
 M^s(S_i)= \binom{n-1}{i}^{-1} \omega_{n-i+s} s! \Phi_{i,s},  \quad i \geq 0. 
\end{displaymath}

The additive kinematic  operator was computed by Schneider \cite{schneider75} and is given by 
\begin{displaymath}
 A(S_i)=\frac{1}{\omega_n} \sum_{k+l=i} \binom{i}{k} S_k \otimes S_l.
\end{displaymath}

Using the diagram above, we get 
\begin{align*}
 a^{\On(n)}_{s_1,s_2}(\Phi_{i,s_1+s_2}) & =\binom{n-1}{i} \frac{1}{\omega_{n-i+s_1+s_2} (s_1+s_2)!} a^{\On(n)}_{s_1,s_2}
\circ M^{s_1+s_2}(S_i)\\
& =\binom{n-1}{i} \frac{1}{\omega_{n-i+s_1+s_2} (s_1+s_2)!} M^{s_1} \otimes M^{s_2} \circ A (S_i)\\ 
& = \binom{s_1+s_2}{s_1}^{-1} \sum_{k+l=i} \frac{(n-k-1)!(n-l-1)!}{(n-i-1)!(n-1)!}\cdot\\
&\qquad \cdot 
\frac{\omega_{n-k+s_1}
\omega_{n-l+s_2}}{\omega_n \omega_{n-i+s_1+s_2}} \Phi_{k,s_1}
\otimes   \Phi_{l,s_2},
\end{align*}
which proves the asserted formula.
\endproof

We remark that we do not obtain the additive kinematic operator on multiples of $Q$ with this method. One possible way to overcome this difficulty is to apply Corollary \ref{cor_relation_kinformulas} to the intersectional kinematic formulas from the next subsection.  Alternatively, one may use Theorem \ref{thm_ftaig} (more precisely the upper square in the diagram) to compute the additive kinematic operator using the knowledge of the convolution algebra from Theorem \ref{mainthm_algebra_tensorvals}. The Poincar\'e pairing, i.e., the vertical arrows in that diagram, may be computed using Corollary \ref{cor_poincare_prod} and \eqref{eq_traces}. We will not go into details.

\subsection{Intersectional kinematic formulas}

In the last few sections, we have introduced the necessary machinery to compute all intersectional kinematic formulas.

Note that for any given ranks $s_1,s_2$, the formula in Theorem \ref{thm_product_tensorcase} is closed. The kinematic
formula $k_{s_1,s_2}^{\On(n)}$ can be explicitly computed with Theorem
\ref{thm_ftaig} and Corollary \ref{cor_poincare_prod}. 

Let us work this out in one of the simplest cases, namely $s_1=s_2=3$, thus showing the last equation in Theorem
\ref{mainthm_kf}. It is clear that there is a formula of the form  
\begin{displaymath}
 k_{3,3}^{\On(n)}(\Phi_{i,6}) = \sum_{k+l=n+i} a_{n,i,k} \Phi_{k,3} \otimes \Phi_{l,3}
\end{displaymath}
with some constants $a_{n,i,k}$ which remain to be determined.

Fix $k,l$ with $k+l=n+i$. Using Corollary
\ref{cor_poincare_prod}, we find 
\begin{align*}
\langle \pd_m^3 \Phi_{k,3},\Phi_{n-k,3}\rangle & =  \frac{1}{72 \pi^3}\binom{n}{k} 
k(n-k)
\frac{\Gamma\left(\frac{k+3}{2}\right)\Gamma\left(\frac{n-k+3}{2}\right)}{\Gamma\left(\frac{n}{2}+1\right)},\\
\langle \pd_m^3 \Phi_{l,3},\Phi_{n-l,3}\rangle & =  \frac{1}{72 \pi^3}\binom{n}{l} 
l(n-l)
\frac{\Gamma\left(\frac{l+3}{2}\right)\Gamma\left(\frac{n-l+3}{2}\right)}{\Gamma\left(\frac{n}{2}+1\right)}
\end{align*}
and therefore 
\begin{multline*}
\langle (\pd_m^3 \otimes \pd_m^3) \circ k_{3,3}^{\On(n)}(\Phi_{i,6}),\Phi_{n-k,3} \otimes \Phi_{n-l,3}\rangle\\
=a_{n,i,k}  \frac{1}{72 \pi^3}\binom{n}{k} 
k(n-k)
\frac{\Gamma\left(\frac{k+3}{2}\right)\Gamma\left(\frac{n-k+3}{2}\right)}{\Gamma\left(\frac{n}{2}+1\right)} 
\frac{1}{72 \pi^3}\binom{n}{l} 
l(n-l)
\frac{\Gamma\left(\frac{l+3}{2}\right)\Gamma\left(\frac{n-l+3}{2}\right)}{\Gamma\left(\frac{n}{2}+1\right)}.
\end{multline*}

On the other hand, by Corollaries \ref{cor_product_tensor_small_rank} and \ref{cor_poincare_prod}, 
\begin{align*}
 &\langle m^* \circ \pd_m^6(\Phi_{i,6}),\Phi_{n-k,3} \otimes \Phi_{n-l,3}\rangle = \langle
\pd_m^6(\Phi_{i,6}),\Phi_{n-k,3} \cdot \Phi_{n-l,3}\rangle \\
&= \frac{(n-k+1)(n-l+1)\Gamma\left(\frac{n-i+1}{2}\right)}{\pi^\frac52
(n-i+4)(n-i+2)(n-i)\Gamma\left(\frac{n-l}{2}\right)\Gamma\left(\frac{n-k}{2}\right)}\cdot\\
&\quad\cdot\left\langle \pd_m^6(\Phi_{i,6}),-32 \Phi_{n-i,6}\pi^3+8Q\Phi_{n-i,4}\pi^2-Q^2
\Phi_{n-i,2}\pi+\frac{1}{12}Q^3\Phi_{n-i,0}\right\rangle\\
&= \frac{1}{207360} \frac{(k-n-1)(i-k-1) \Gamma\left(\frac{n+1}{2}\right)
(i+1)(i-1)(i-3)}{\pi^5 \Gamma\left(\frac{i+1}{2}\right) \Gamma\left(\frac{n-k}{2}\right)
\Gamma\left(\frac{k-i}{2}\right)}.
\end{align*}
From this, the explicit value of $a_{n,i,k}$ given in the theorem follows. 

Let us now prove the formula for $k_{3,2}^{\On(n)}$. There are constants $a_{n,i,k},b_{n,i,k}$ such that 
\begin{displaymath}
 k_{3,2}^{\On(n)}(\Phi_{i,5}) = \sum_{k+l=n+i} \Phi_{k,3} \otimes \left(a_{n,i,k}\Phi_{l,2}+b_{n,i,k} Q
\Phi_{l,0}\right).
\end{displaymath}
Fix $k,l$ with $k+l=n+i$. By Corollary \ref{cor_poincare_prod} and the subsequent remark, we can compute the relevant
Poincar\'e pairings. For instance 

\begin{align*}
\langle \pd_m^5 \Phi_i^5,Q\Phi_{n-i,3}\rangle & = \langle \pd_m^3 \tr \Phi_i^5,\Phi_{n-i,3}\rangle \\
& = \frac{n-i+3}{40 \pi} \langle \pd_m^3 \Phi_i^3,\Phi_{n-i,3}\rangle \\
& = \frac{n-i+3}{40 \pi}  \frac{1}{72 \pi^3 }\binom{n}{i} i(n-i)
\frac{\Gamma\left(\frac{i+3}{2}\right)\Gamma\left(\frac{n-i+3}{2}\right)}{\Gamma\left(\frac{n}{2}+1\right)}.
\end{align*}

The products $\Phi_{n-k,3} \cdot \Phi_{n-l,2}$ and $\Phi_{n-k,3} \cdot Q\Phi_{n-l,0}=Q (\Phi_{n-k,3}\cdot \mu_{n-l})$
were computed in Corollary \ref{cor_product_tensor_small_rank} and in Subsection
\ref{subsec_croftonformulas} respectively. 

From the two equations 
\begin{align*}
 \langle (\pd_m^3 \otimes \pd_m^2) \circ k_{3,2}^{\On(n)}(\Phi_{i,5}),\Phi_{n-k,3} \otimes \Phi_{n-l,2}\rangle & =
\langle \pd_m^5 \Phi_{i,5},\Phi_{n-k,3} \cdot \Phi_{n-l,2}\rangle,\\
 \langle (\pd_m^3 \otimes \pd_m^2) \circ k_{3,2}^{\On(n)}(\Phi_{i,5}),\Phi_{n-k,3} \otimes Q\Phi_{n-l,0}\rangle & =
\langle \pd_m^5 \Phi_{i,5},\Phi_{n-k,3} \cdot Q\Phi_{n-l,0}\rangle
\end{align*}
we obtain a linear system of equations for the unknowns $a_{n,k,i}, b_{n,k,i}$ whose unique solution is 
\begin{align*}
 a_{n,k,i} & = \frac{(i+1)\Gamma\left(\frac{l+1}{2}\right)\Gamma\left(\frac{k}{2}\right)}{40 \pi
(k+1)l \Gamma\left(\frac{n+1}{2}\right) \Gamma\left(\frac{i}{2}\right)} 4 \pi (i-3)\\
 b_{n,k,i} & = \frac{(i+1)\Gamma\left(\frac{l+1}{2}\right)\Gamma\left(\frac{k}{2}\right)}{40 \pi
(k+1)l \Gamma\left(\frac{n+1}{2}\right) \Gamma\left(\frac{i}{2}\right)} (n-k+3).
\end{align*}
This finishes the proof of the second equation in Theorem \ref{mainthm_kf}. The first equation is obtained in a
similar way, with a linear system of four equations to solve.  


\end{document}